\documentclass[11pt]{amsart}
\usepackage{latexsym}
\usepackage{amsmath}
\usepackage{amssymb}

\newcommand{\eproof}{\mbox{\ }\hfill $\Box$ \par \vskip 10pt}

\newtheorem{Theorem}{Theorem}[section]
\newtheorem{lemma}[Theorem]{Lemma}
\newtheorem{prop}[Theorem]{Proposition}
\newtheorem{rem}[Theorem]{Remark}

\numberwithin{equation}{section}

\def\R{{\mathbb R}}

\def\cal{\mathcal}

\begin{document}

\title[Semiclassical resolvent estimates]{Semiclassical resolvent estimates for the magnetic Schr\"odinger operator}

\author[G. Vodev]{Georgi Vodev}

\address {Universit\'e de Nantes, Laboratoire de Math\'ematiques Jean Leray, 2 rue de la Houssini\`ere, BP 92208, 44322 Nantes Cedex 03, France}
\email{Georgi.Vodev@univ-nantes.fr}

\date{}

\begin{abstract} We obtain semiclassical resolvent estimates for the
Schr\"odinger operator $(ih\nabla +b)^2+V$ in $\mathbb{R}^d$, $d\ge 3$, where $0<h\ll 1$ is a semiclassical parameter, 
$V$ and $b$ are real-valued electric and 
magnetic potentials independent of $h$. 
If $V\in L^\infty(\mathbb{R}^d)$, $b\in L^\infty(\mathbb{R}^d;\mathbb{R}^d)$, ${\rm div}\,b\in L^\infty(\mathbb{R}^d)$ 
satisfy $V(x)=\mathcal{O}\left(|x|^{-1-\epsilon}\right)$, $b(x)=\mathcal{O}\left(|x|^{-1-\epsilon}\right)$,
${\rm div}\,b(x)=\mathcal{O}\left(|x|^{-1-\epsilon}\right)$, $\epsilon>0$, for $|x|\gg 1$, 
we prove that the norm of the weighted resolvent is bounded by $\exp\left(Ch^{-2}\log(h^{-1})\right)$, $C>0$. 
 We get better
resolvent bounds for electric potentials which are H\"older with respect to the radial variable and 
magnetic potentials which are H\"older with respect to the space variable. 
For long-range electric potentials which are Lipschitz with respect to the radial variable and long-range 
magnetic potentials which are Lipschitz with respect to the space variable we obtain a resolvent bound of
the form $\exp\left(Ch^{-1}\right)$, $C>0$.

\quad

\end{abstract} 

\maketitle

\setcounter{section}{0}
\section{Introduction and statement of results}

In this paper we  
study the resolvent of the Schr\"odinger operator
\begin{equation*}
P(h)=(ih\nabla +b(x))^2+V(x)\quad \mbox{in} \quad\mathbb{R}^d,\, d\ge 3,
\end{equation*}
where $0<h\ll 1$ is a semiclassical parameter, $\nabla$ is the gradient, $V\in L^\infty(\mathbb{R}^d;\mathbb{R})$ and 
$b=(b_1,...,b_d)\in L^\infty(\mathbb{R}^d;\mathbb{R}^d)$ are electric and
magnetic potentials, respectively. 
We define the self-adjoint realization of the operator $(ih\nabla+b)^2+V$ on the Hilbert space $L^2(\mathbb{R}^d)$, 
which again will be denoted by $P(h)$, via the quadratic form $q$ on the Sobolev space $H^1(\R^d)$ 
defined by
\begin{equation*}
q (u,v)= \int_{\mathbb{R}^d}((ih\nabla+b)u\cdot\overline{(ih\nabla+b)v}+Vu\overline{v})dx.
\end{equation*}
Then the domain of $P(h)$ satisfies $D(P(h))\subset H^1(\R^d)$. We refer to Appendix A of both papers 
\cite{kn:LSV} and \cite{kn:MPS} for more details.

We are interested in bounding the quantity
\begin{equation*}
g_s^\pm(h,\varepsilon):=\log\left\|\langle x\rangle^{-s}(P(h)-E\pm i\varepsilon)^{-1}\langle x\rangle^{-s}
\right\|_{L^2\to L^2}
\end{equation*}
from above by an explicit function of $h$, independent of $\varepsilon$. 
Here $\langle x\rangle:=(|x|^2+1)^{1/2}$, 
$L^2:=L^2(\mathbb{R}^d)$, $0<\varepsilon<1$, $s>1/2$ is independent of $h$ and $E>0$ is a fixed energy level independent of $h$.
Our goal is to obtain the best possible upper bounds for $g_s^\pm(h,\varepsilon)$ for potentials which are Lipschitz, H\"older
or just $L^\infty$. There have been recently many papers studying this problem when $b\equiv 0$. To our best knoweledge, no such
results exist for non-trivial magnetic potentials when $d\ge 2$ and it seems that this paper is the first one where upper bounds for
$g_s^\pm$ are proved in this case.  When $d=1$ and $b$ not identically zero, 
semiclassical resolvent bounds have been recently proved in \cite{kn:LS1} for a very large class of electric and magnetic potentials. 
Note also that sharp high-frequency resolvent bounds for the operator $P(1)$ with $L^\infty$ potentials are proved in 
\cite{kn:MPS} when $d\ge 2$ and in \cite{kn:LSV} and \cite{kn:V1} when $d\ge 3$. In \cite{kn:G} exponential high-frequency 
resolvent bounds for the operator $P(1)$ with smooth potentials on non-compact Riemannian manifolds have been recently proved,
extending the results in \cite{kn:CV}. 

We will be looking for bounding $g_s^\pm$ for the largest possible class of electric and magnetic 
potentials that our method allows to cover.
To describe it we introduce the polar coordinates $r=|x|$, $w=x/|x|$. We suppose that $V=V_L+V_S$, $b=b^L+b^S$, where 
$V_L, b^L$ (resp. $V_S, b^S$) are long-range (resp. short-range) electric and magnetic potentials satisfying the following conditions.
We suppose that $V_L$ satisfies
\begin{equation}\label{eq:1.1}
V_L(rw)\le p(r),
\end{equation}
where $p>0$ is a decreasing function such that $p(r)\to 0$ as $r\to\infty$. 
We also suppose that $V_L(rw)$ is Lipschitz with respect to the radial variable $r$ and its first derivative
satisfies the bound
\begin{equation}\label{eq:1.2}
\partial_rV_L(rw)\le C(r+1)^{-1-\rho}.
\end{equation}
Hereafter $C>0$ and $0<\rho\ll 1$ denote positive constants that may change from line to line. 
The short-range part $V_S\in L^\infty(\mathbb{R}^d)$ is supposed to satisfy 
\begin{equation}\label{eq:1.3}
|V_S(rw)|\le C(r+1)^{-1-\rho}.
\end{equation}
The long-range part $b^L\in L^\infty(\mathbb{R}^d;\mathbb{R}^d)$
is supposed to be Lipschitz with respect to $r$ and satisfies 
\begin{equation}\label{eq:1.4}
|\partial_r^kb^L(rw)|\le C(r+1)^{-k-\rho},\quad k=0,1.
\end{equation}
The short-range part $b^S\in L^\infty(\mathbb{R}^d;\mathbb{R}^d)$ satisfies  
\begin{equation}\label{eq:1.5}
|b^S(rw)|\le C(r+1)^{-1-\rho}.
\end{equation}
Finally, we suppose that ${\rm div}\,b^S\in L^\infty(\mathbb{R}^d)$ exists and satisfies
\begin{equation}\label{eq:1.6}
|{\rm div}\,b^S(rw)|\le C(r+1)^{-1-\rho}.
\end{equation}
This condition is fulfilled if each $b_j^S$ is Lipschitz with respect to the variable $x_j$
and $\partial_{x_j}b_j^S(x)=\mathcal{O}\left(|x|^{-1-\rho}\right)$ for $|x|\gg 1$. 

Throughout this paper, given two vectors $a,b\in \mathbb{R}^d$, $a\cdot b$ will denote the scalar product. 
Our first result is the following

\begin{Theorem} \label{1.1}
In addition to the above conditions we assume either the condition
\begin{equation}\label{eq:1.7}
|w\cdot b^L(rw)|\le C(r+1)^{-1-\rho},
\end{equation}
or we suppose that the function $x\cdot b^L(x)$ is Lipschitz in $x$ and satisfies  
\begin{equation}\label{eq:1.8}
|\partial_x^\beta(x\cdot b^L(x))|\le C,\quad |\beta|=1.
\end{equation}
Then, there exist constants $C>0$, $0<h_0\ll 1$, independent of $h$ and $\varepsilon$, such 
that the bound
\begin{equation}\label{eq:1.9}
g_s^\pm(h,\varepsilon)\le Ch^{-2}\log(h^{-1})
\end{equation}
holds for all $0<h\le h_0$. If $V_S\equiv 0$ and $b^S\equiv 0$, 
under the condition (\ref{eq:1.8}), we have the better bound
\begin{equation}\label{eq:1.10}
g_s^\pm(h,\varepsilon)\le Ch^{-1}.
\end{equation}
\end{Theorem}

Note that 
the bound (\ref{eq:1.10}) is proved in \cite{kn:LS1} when $d=1$, while for $d\ge 3$ it seems to be new. 
Proving (\ref{eq:1.10}) when $d=2$ in the presence of a magnetic potential remains an open problem. 
When $b\equiv 0$, however, the bound (\ref{eq:1.10}) is well-known. Indeed, it was proved in \cite{kn:D} 
when $d\ge 3$ for potentials satisfying (\ref{eq:1.1})
with $p(r)=C(r+1)^{-\rho}$, and extended to the general case in \cite{kn:GS1} and \cite{kn:V3}. 
When $d=2$ the bound (\ref{eq:1.10}) is proved in \cite{kn:S1}
for potentials $V\in C^1(\mathbb{R}^2)$ satisfying (\ref{eq:1.1})
with $p(r)=C(r+1)^{-\rho}$ as well as the condition
\begin{equation}\label{eq:1.11}
|\nabla V(x)|\le C\langle x\rangle^{-1-\rho}.
\end{equation}
The bound (\ref{eq:1.10}) has been recently proved in \cite{kn:O} for long-range Lipschitz potentials having singularities at the origin
in all dimensions $d\ge 2$. 
When $d=1$ the bound (\ref{eq:1.10}) is proved in \cite{kn:DS} for potentials $V\in L^1(\mathbb{R})$ and in \cite{kn:LS1} for more general
measure potentials. 
It turns out that (\ref{eq:1.10}) also holds for $V\in L^\infty(\mathbb{R}^d)$, $d\ge 2$, provided $V$ is compactly supported 
and depends only on the radial variable $r$. Indeed, this is proved in \cite{kn:DGS} when $d\ge 2$.
A simpler proof of this result is given in \cite{kn:V6} when $d\ge 3$. Without a radial symmetry or some regularity
of the potential, however, the bound (\ref{eq:1.10}) seems to fail even for compactly supported potentials. 
Indeed, for potentials $V\in L^\infty$ satisfying  
\begin{equation}\label{eq:1.12}
|V(x)|\le C\langle x\rangle^{-\delta}
\end{equation}
where $\delta>1$, the following bounds are known to hold:
\begin{equation}\label{eq:1.13}
g_s^\pm(h,\varepsilon)\le Ch^{-4/3}\log(h^{-1})\quad\mbox{if}\quad \delta>2,
\end{equation}
\begin{equation}\label{eq:1.14}
g_s^\pm(h,\varepsilon)\le C_\epsilon h^{-\frac{2\delta}{2\delta-1}-\epsilon}\quad\mbox{if}\quad 1<\delta\le 2,
\end{equation}
for every $0<\epsilon\ll 1$. When $d\ge 3$ the bound (\ref{eq:1.13}) is proved in \cite{kn:GS2} and in \cite{kn:V2} if $\delta>3$. 
Previously, (\ref{eq:1.13}) was proved for compactly supported $L^\infty$ potentials 
in \cite{kn:KV} and \cite{kn:S2}. The bound (\ref{eq:1.13}) is extended in \cite{kn:V4} to more general Riemannian manifolds. 
The bounds (\ref{eq:1.13}) and (\ref{eq:1.14}) have been recently proved in 
\cite{kn:S3} when $d\ge 2$ for potentials $V=V_L+V_S$ with $V_L$ satisfying (\ref{eq:1.1}) and (\ref{eq:1.2}),
and $V_S$ satisfying (\ref{eq:1.12}). In fact, the conditions on the behavior of $V_L$ and $V_S$ at infinity in this paper 
are a little bit more general than this. Moreover, the potential in \cite{kn:S3} is allowed to have singularities at the origin.
It turns out that (\ref{eq:1.13}) and (\ref{eq:1.14})
can be improved for potentials $V\in L^\infty$ depending only on $r$. Indeed, for radial potentials
satisfying (\ref{eq:1.12}) it is shown in \cite{kn:V5} and \cite{kn:V6} that when $d\ge 3$ we have the bounds
\begin{equation}\label{eq:1.15}
g_s^\pm(h,\varepsilon)\le Ch^{-\frac{\delta}{\delta-1}}\quad\mbox{if}\quad \delta>4,
\end{equation}
\begin{equation}\label{eq:1.16}
g_s^\pm(h,\varepsilon)\le Ch^{-4/3}\quad\mbox{if}\quad 2<\delta\le 4,
\end{equation}
\begin{equation}\label{eq:1.17}
g_s^\pm(h,\varepsilon)\le Ch^{-\frac{2\delta}{2\delta-1}}
\left(\log(h^{-1})\right)^{\frac{\delta+1}{2\delta-1}}\quad\mbox{if}\quad 1<\delta\le 2.
\end{equation}
The above bounds can be improved if some small regularity of
the potential is assumed. 
To desribe these results, given $0<\alpha\le 1$ and $k>0$, 
we introduce the space $C_k^\alpha(\overline{\mathbb{R}^+})$, $\mathbb{R}^+:=(0,\infty)$, of all H\"older functions $a$ such that
 \begin{equation*}
 \sup_{r'\ge 0:\,0<|r-r'|\le 1}\frac{|a(r)-a(r')|}{|r-r'|^\alpha}\le C(r+1)^{-k},\quad\forall r\in \overline{\mathbb{R}^+}.
 \end{equation*}
We suppose that the function $V(r,w):=V(rw)$ satisfies the condition
\begin{equation}\label{eq:1.18}
V(\cdot,w)\in C_{3+\rho}^\alpha(\overline{\mathbb{R}^+}),\quad 0<\alpha<1,
\end{equation}
uniformly in $w\in \mathbb{S}^{d-1}$. For potentials $V$ satisfying (\ref{eq:1.1}) and (\ref{eq:1.18}) and $d\ge 3$ the following 
 bound is known to hold:
\begin{equation}\label{eq:1.19}
g_s^\pm(h,\varepsilon)\le Ch^{-4/(\alpha+3)}\log(h^{-1}).
\end{equation}
The bound (\ref{eq:1.19}) is proved in \cite{kn:GS2} and in \cite{kn:V3} if $\rho=1$. It has been 
extended in \cite{kn:V3} to the case $d=2$ as well as to exterior domains for potentials which are 
H\"older with respect to the space variable $x$. 
The bound (\ref{eq:1.19}) is improved in
\cite{kn:V5} for radial potentials satisfying (\ref{eq:1.1}) and (\ref{eq:1.18}) with $\rho=0$ to 
\begin{equation}\label{eq:1.20}
g_s^\pm(h,\varepsilon)\le Ch^{-4/(\alpha+3)}.
\end{equation}
The proof of the above bounds is based on global Carleman estimates with suitable phase and weight functions. 
In order that the Carleman estimates work the phase and weight functions must satisfy some conditions.
Therefore, the task consists of constructing explicitly these functions in such a way that we get the best possible
bound for $g_s^\pm$. In fact, there are many ways to do so and the class the potential belongs to is determined by the way one does it.

To prove Theorem 1.1 we follow the same strategy as in the case $b\equiv 0$. The Carleman estimates for the magnetic Schr\"odinger
operator, however, are much more difficult to handle with. The main difference is that conjugating this operator with
$e^{\varphi/h}$, $\varphi>0$ being a phase function depending only on $r$, makes appear a complex-valued 
effective potential of the form 
\begin{equation*}
-2i\varphi'(r)w\cdot b(rw),
\end{equation*}
where $\varphi'>0$ denotes the first derivative of $\varphi$. 
This term is not easy to handle with even in the case when $b^S\equiv 0$.
One way to treat it is to consider it as a short-range potential and to add it to $V_S$. To do so, however, 
we need the condition (\ref{eq:1.7}). It allows us to make the Carleman estimates work. However, this approach 
does not lead to better bounds
than what we have in (\ref{eq:1.9}). Therefore, we propose another way that requires the condition (\ref{eq:1.8}). 
This approach allows us to improve the bound (\ref{eq:1.9}) if $b_S\equiv 0$ or $b_S$ is H\"older. To be more precise, 
 given $0<\alpha\le 1$ and $k>0$, 
we introduce the space $C_k^\alpha(\mathbb{R}^d)$ of all H\"older functions $a$ such that
 \begin{equation*}
 \sup_{x'\in\mathbb{R}^d:\,0<|x-x'|\le 1}\frac{|a(x)-a(x')|}{|x-x'|^\alpha}\le C\langle x\rangle^{-k},\quad\forall x\in\mathbb{R}^d.
 \end{equation*}
  Note that the case $\alpha=1$ corresponds to the Lipschitz functions. We suppose that
  \begin{equation}\label{eq:1.21}
V_S(\cdot,w)\in C_{3+\rho}^\alpha(\overline{\mathbb{R}^+}),\quad 0<\alpha\le 1,\quad \rho>0,
\end{equation}
uniformly in $w$, and
\begin{equation}\label{eq:1.22}
b^S\in C_{3/2+\rho_1}^{\alpha'}(\mathbb{R}^d;\mathbb{R}^d),\quad 0<\alpha'\le 1,\quad \rho_1>0.
\end{equation}
 We have the following

\begin{Theorem} \label{1.2}
Assume that $V_L$ satisfies (\ref{eq:1.1}), (\ref{eq:1.2}), $b^L$ satisfies (\ref{eq:1.4}), (\ref{eq:1.8}), $V_S$ satisfies 
(\ref{eq:1.21}) and (\ref{eq:1.1}) with possibly a new decreasing function $p$ tending to zero, $b^S$ satisfies (\ref{eq:1.6}), (\ref{eq:1.22}).
We also suppose that $b^L$ and $b^S$ satisfy
\begin{equation}\label{eq:1.23}
|b^L(rw)|+|b^S(rw)|\le C(r+1)^{-\rho_2},
\end{equation}
with some $\rho_2>0$ such that $\rho_1+\rho_2>3/2$.
Then we have the bound
\begin{equation}\label{eq:1.24}
g_s^\pm(h,\varepsilon)\le Ch^{-n}\log(h^{-1}),
\end{equation}
where
\begin{equation*}
n=\max\left\{\frac{4}{\alpha+3},\frac{2}{\alpha'+1}\right\}.
\end{equation*}
\end{Theorem}

When $\alpha'\le\frac{1}{2}$ the bound (\ref{eq:1.24}) also holds for $L^\infty$ potentials $V_S$. More precisely, we have the following

\begin{Theorem} \label{1.3}
Assume that $V_L$ satisfies (\ref{eq:1.1}), (\ref{eq:1.2}), $b^L$ satisfies (\ref{eq:1.4}), (\ref{eq:1.8})
and (\ref{eq:1.23}), $b^S$ satisfies (\ref{eq:1.6}), 
(\ref{eq:1.22}) 
and (\ref{eq:1.23}), $V_S$ satisfies 
\begin{equation}\label{eq:1.25}
|V_S(rw)|\le C(r+1)^{-\delta},
\end{equation}
with $\delta>1$. 
Then we have the bounds
\begin{equation}\label{eq:1.26}
g_s^\pm(h,\varepsilon)\le Ch^{-2/(\alpha'+1)}\log(h^{-1})\quad\mbox{if}\quad \delta>\frac{1}{1-\alpha'},\,\,\alpha'\le\frac{1}{2},
\end{equation}
\begin{equation}\label{eq:1.27}
g_s^\pm(h,\varepsilon)\le C_\epsilon h^{-\frac{2\delta}{2\delta-1}-\epsilon}\quad\mbox{if}\quad \alpha'\ge 
\frac{\delta-1}{\delta},\,\,1<\delta\le 2,
\end{equation}
for every $0<\epsilon\ll 1$ independent of $h$.
\end{Theorem}

The paper is organized as follows. In Section 2 we construct the phase and weight functions needed for the 
proof of the Carleman estimates. In Section 3 we prove Carleman estimates for the magnetic Schr\"odinger operator under the conditions
of Theorem \ref{1.1} needed for the proof of the bound (\ref{eq:1.9}). In Section 4 we improve the Carleman estimates
when $V_S\equiv 0$ and $b^S\equiv 0$, which leads to the bound (\ref{eq:1.10}). 
In Section 5 we adapt the Carleman estimates to H\"older potentials $V_S$ and $b^S$, which leads to the bound (\ref{eq:1.24}).
In Section 6 we adapt the Carleman estimates to $L^\infty$ potentials $V_S$ and 
H\"older potentials $b^S$, which leads to the bounds (\ref{eq:1.26}) and (\ref{eq:1.27}).
 In Section 7 we show that the Carleman estimates imply the resolvent ones.
Finally, Section 8 is devoted to the proof of Lemma \ref{7.1}. 

\section{Construction of the phase and weight functions} 
We first construct the weight function $\mu$ as follows:
\begin{equation*}
\mu(r)=
\left\{
\begin{array}{lll}
 r^2+r&\mbox{for}& 0\le r\le 1,\\
  r^2+r^{2-\ell}&\mbox{for}& 1\le r\le a,\\
 (a^2+a^{2-\ell})(1+(a+1)^{-2s+1}-(r+1)^{-2s+1})&\mbox{for}& r\ge a,
\end{array}
\right.
\end{equation*}
where $a\gg 1$ is either independent of $h$ or it is of the form $a=h^{-m}$, $m>0$.
The parameters $\ell$ and $s$ are  independent of $h$ such that $0<\ell\ll 1$, $1/2<s<1$. 
Clearly, the first derivative of $\mu$ is given by
\begin{equation*}
\mu'(r)=
\left\{
\begin{array}{lll}
 2r+1&\mbox{for}& 0\le r<1,\\
 2r+(2-\ell)r^{1-\ell}&\mbox{for}& 1<r<a,\\
 (2s-1)(a^2+a^{2-\ell})(r+1)^{-2s}&\mbox{for}& r>a.
\end{array}
\right.
\end{equation*}
We have the following

\begin{lemma} \label{2.1}
We have the identity
\begin{equation}\label{eq:2.1}
2r^{-1}\mu(r)-\mu'(r)=
\left\{
\begin{array}{lll}
 1&\mbox{for}& 0<r<1,\\
 \ell r^{1-\ell}&\mbox{for}& 1<r<a.
\end{array}
\right.
\end{equation}
For all $r>a$, we have the lower bound
\begin{equation}\label{eq:2.2}
2r^{-1}\mu(r)-\mu'(r)\ge a^2(r+1)^{-1}.
\end{equation}
For all $r\neq 1$, $r\neq a$, we have
\begin{equation}\label{eq:2.3}
\mu'(r)\ge (r+1)^{-2s}.
\end{equation}
\end{lemma}

{\it Proof.} The identity (\ref{eq:2.1}) is obvious, while (\ref{eq:2.2}) follows from the inequalities
\begin{equation*}
2r^{-1}\mu(r)\ge 2(a^2+a^{2-\ell})(r+1)^{-1},\quad \mu'(r)\le (2s-1)(a^2+a^{2-\ell})(r+1)^{-1}.
\end{equation*}
It is also clear that the inequality (\ref{eq:2.3}) holds as long as $a^2\ge (2s-1)^{-1}$. 
\eproof

We will now construct a phase function
$\psi\in C^1([0,+\infty))$ such that $\psi(0)=0$ and $\psi(r)>0$ for $r>0$. 
We define the first derivative of $\psi$ by
\begin{equation*}
\psi'(r)=
\left\{
\begin{array}{lll}
 c_\nu&\mbox{for}& 0\le r\le 1,\\
 r^{-1}(1+r^{-\ell})^{-1/2}(1-(r+1)^{-\nu})&\mbox{for}& 1\le r\le a,\\
 K_a(r+1)^{-1-\kappa}&\mbox{for}& r\ge a,
\end{array}
\right.
\end{equation*}
where $0<\nu,\kappa\ll 1$ are independent of $h$ to be fixed later on, and
\begin{equation*}
c_\nu=2^{-1/2}(1-2^{-\nu}),\quad K_a=a^{-1}(a+1)^{1+\kappa}(1+a^{-\ell})^{-1/2}(1-(a+1)^{-\nu})=\mathcal{O}\left(a^{\kappa}\right),
\end{equation*}
are chosen in such a way that the function $\psi'$ is Lipschitz. 
It is easy to check that the first derivative of $\psi'$ satisfies
\begin{equation*}
\psi''(r)=
\left\{
\begin{array}{lll}
 0&\mbox{for}& 0\le r<1,\\
 \mathcal{O}(r^{-2})&\mbox{for}& 1<r<a,\\
 \mathcal{O}(a^{\kappa}r^{-2-\kappa})&\mbox{for}& r>a.
\end{array}
\right.
\end{equation*}

\begin{lemma} \label{2.2}
For all $r\ge 0$ we have the bound
\begin{equation}\label{eq:2.4}
\psi(r)\lesssim \log a.
\end{equation}
\end{lemma}

{\it Proof.} We have
\begin{equation*}
\begin{split}
\max\psi &=\int_0^1\psi'(r)dr+\int_1^a\psi'(r)dr+\int_a^\infty\psi'(r)dr\\
&\le c_\nu\int_0^1dr+\int_1^ar^{-1}dr+K_a\int_a^\infty (r+1)^{-1-\kappa}dr\\
&=c_\nu+\log a+\kappa^{-1}(a+1)^{-\kappa}K_a,
\end{split}
\end{equation*}
which confirms (\ref{eq:2.4}).
\eproof
The next lemma will play a crucial role in the proof of the Carleman estimates. 

\begin{lemma} \label{2.3}
For $0<r<a$, $r\neq 1$, we have the inequality
\begin{equation}\label{eq:2.5}
\left(\mu\psi'^2\right)'(r)\ge C_\nu(r+1)^{-1-\nu},
\end{equation}
with some constant $C_\nu>0$. 
For all $r>a$ we have the inequality
\begin{equation}\label{eq:2.6}
\left(\mu\psi'^2\right)'(r)\ge -\widetilde C_\kappa a^{-3+2s}\mu'(r), 
\end{equation}
with some constant $\widetilde C_\kappa>0$. 
\end{lemma}

{\it Proof.} For $0<r<1$ we have 
 \begin{equation*}
 \left(\mu\psi'^2\right)'(r)=c_\nu^2(2r+1),
 \end{equation*}
 which clearly implies (\ref{eq:2.5}) in this case. For $1<r<a$ we have
 \begin{equation*}
 \begin{split}
 \left(\mu\psi'^2\right)'(r)&=((1-(r+1)^{-\nu})^2)'\\
 &=2\nu(1-(r+1)^{-\nu})(r+1)^{-1-\nu}\\
 &\ge 2\nu(1-2^{-\nu})(r+1)^{-1-\nu}.
 \end{split}
 \end{equation*}
For $r>a$, we have $\mu(r)\lesssim (r+1)^{2s}\mu'(r)$. Hence
 \begin{equation*}
 \begin{split}
\left(\mu\psi'^2\right)'(r)&=\mu'\psi'^2+2\mu\psi'\psi''\\
&\ge -C(r+1)^{2s}\mu'\psi'|\psi''|\\
&\ge -Ca^{2\kappa}(r+1)^{-3-2\kappa+2s}\mu'\\
&\ge -Ca^{-3+2s}\mu',
\end{split}
 \end{equation*}
with some constant $C>0$, which confirms (\ref{eq:2.6}).
\eproof

\section{Carleman estimates in the general case}

Let $\mu$ and $\psi$ be the functions introduced in Section 2 and set $\varphi=\tau\psi$, where
$\tau=\tau_0h^{-1}$, 
$\tau_0\gg 1$ being a parameter independent of $h$ to be chosen later on. 
In this section we prove the following

\begin{Theorem} \label{3.1}
Under the conditions of Theorem \ref{1.1}, there are constants $C>0$ and $0<h_0\ll 1$ such that for all $0<h\le h_0$, 
$0<\varepsilon\le 1$, $s>1/2$, and for all functions
$f\in D(P(h))$ such that $(\zeta\partial_r)^kf\in H^1(\mathbb{R}^d)$, $k=0,1$, and satisfying 
\begin{equation*}
\langle x\rangle^{s}(P(h)-E\pm i\varepsilon)f\in L^2(\mathbb{R}^d),
\end{equation*}
 we have the estimate 
 \begin{equation}\label{eq:3.1}
 \begin{split}
\|\langle x\rangle^{-s}e^{\varphi/h}f\|_{L^2}&\le Cah^{-1}\|\langle x\rangle^{s}e^{\varphi/h}(P(h)-E\pm i\varepsilon)f\|_{L^2}\\ 
&+Ca\left(\tau h^{-1}\varepsilon\right)^{1/2}\|e^{\varphi/h}f\|_{L^2}\\
&+Ca\left(\tau^{-1}h^{-1}\varepsilon \right)^{1/2}\|e^{\varphi/h}{\cal D}_rf\|_{L^2},
\end{split}
\end{equation}
where we have set ${\cal D}_r=ih\partial_r$ and $\zeta(r)=r(1+r^2)^{-1/2}$.
\end{Theorem}

{\it Proof.} Clearly, it suffices to prove (\ref{eq:3.1}) for $0<s-1/2\ll 1$ as this would imply (\ref{eq:3.1}) for
all $s>1/2$. So, in what follows $s$ will be as in Section 2 with additional restrictions made later on. 
We write the operator $P(h)$ in the form
\begin{equation*}
P(h)=-h^2\Delta+ihb^L\cdot\nabla+ih\nabla\cdot b^L+2ihb^S\cdot\nabla+\widetilde V_L+\widetilde V_S,
\end{equation*}
where
\begin{align*}
&\widetilde V_L=|b^L|^2+V_L,\\
&\widetilde V_S=ih{\rm div}\,b^S+2b^S\cdot b^L+|b^S|^2+V_S.
\end{align*}
Note that our conditions guarantee that 
\begin{equation*}
|\widetilde V_S(x)|\lesssim\langle x\rangle^{-1-\rho},
\quad \partial_r\widetilde V_L(x)\lesssim\langle x\rangle^{-1-\rho}.
\end{equation*}
Moreover, $\widetilde V_L$ satisfies (\ref{eq:1.1}) with a new decreasing function $p$. 

We will write $P(h)$ in the polar coordinates $(r,w)\in\mathbb{R}^+\times\mathbb{S}^{d-1}$,  
$r=|x|$, $w=x/|x|$. Recall that  $L^2(\mathbb{R}^d)=L^2(\mathbb{R}^+\times\mathbb{S}^{d-1}, r^{d-1}drdw)$. 
In what follows in this section we denote by $\|\cdot\|$ and $\langle\cdot,\cdot\rangle$
the norm and the scalar product in $L^2(\mathbb{S}^{d-1})$. We will make use of the identity
\begin{equation}\label{eq:3.2}
 r^{(d-1)/2}\Delta  r^{-(d-1)/2}=\partial_r^2+r^{-2}\Delta_w-(d-1)(d-3)(2r)^{-2}
\end{equation}
where $\Delta_w$ denotes the negative Laplace-Beltrami operator
on $\mathbb{S}^{d-1}$. We also need the following lemma the proof of which is given in Appendix A.

\begin{lemma} \label{3.2} 
We have the formula
\begin{equation}\label{eq:3.3}
 r^{(d-1)/2}\partial_{x_j} r^{-(d-1)/2}=w_j\partial_r+r^{-1}q_j(w,\partial_w)
\end{equation}
where $q_j$ is a first-order differential operator on $\mathbb{S}^{d-1}$ which is antisymmetric with respect to the scalar product
in $L^2(\mathbb{S}^{d-1})$. 
\end{lemma}

Clearly, we have the estimates
\begin{equation}\label{eq:3.4}
 \|q_j(w,\partial_w)v\|\lesssim \|(-\Delta_w)^{1/2}v\|+\|v\|\quad\mbox{for all}\quad v\in H^1(\mathbb{S}^{d-1}).
\end{equation}
Set $\Lambda=-h^2\Delta_w$, $\mathcal{D}_r=ih\partial_r$ and $Q_j=ihq_j(w,\partial_w)$. Then the operator $Q_j$
is symmetric with respect to the scalar product
in $L^2(\mathbb{S}^{d-1})$. Moreover, (\ref{eq:3.4}) implies the estimate
\begin{equation}\label{eq:3.5}
 \|Q_jv\|\lesssim \|\Lambda^{1/2}v\|+h\|v\|\quad\mbox{for all}\quad v\in H^1(\mathbb{S}^{d-1}).
\end{equation}
Set 
\begin{align*}
&\mathcal{P}(h):=r^{(d-1)/2}(P(h)-E)r^{-(d-1)/2},\\
&\mathcal{P}_\varphi(h):=e^{\varphi/h}\mathcal{P}(h)e^{-\varphi/h}.
\end{align*}
Using (\ref{eq:3.2}) and (\ref{eq:3.3}) we can write the operator ${\cal P}(h)$ in the coordinates $(r,w)$ as follows
\begin{equation*}
\begin{split}
\mathcal{P}(h)&=\mathcal{D}_r^2+r^{-2}\Lambda-E+h^2(d-1)(d-3)(2r)^{-2}+\widetilde V_L+\widetilde V_S\\
&+\sum_{j=1}^dw_j(b_j^L\mathcal{D}_r+\mathcal{D}_rb_j^L)+2\sum_{j=1}^dw_jb_j^S\mathcal{D}_r\\
&+r^{-1}\sum_{j=1}^d\left(b_j^LQ_j+Q_jb_j^L\right)+2r^{-1}\sum_{j=1}^db_j^SQ_j.
\end{split}
\end{equation*}
We now define the function $\sigma$ as follows. We put $\sigma=0$ if (\ref{eq:1.7}) holds and we put 
\begin{equation*}
\sigma(r,w)=\sum_{j=1}^dw_jb_j^L(rw)
\end{equation*}
if (\ref{eq:1.8}) holds. 
Now we can write the operator ${\cal P}(h)$ as follows
\begin{equation*}
\begin{split}
\mathcal{P}(h)&=\mathcal{D}_r^2+r^{-2}\Lambda-E+h^2(d-1)(d-3)(2r)^{-2}+\widetilde V_L+\widetilde V_S\\
&+2\sum_{j=1}^dw_jb_j\mathcal{D}_r+ih\sum_{j=1}^dw_j\partial_r(b_j^L)\\
&+r^{-1}\sum_{j=1}^d\left(b_j^LQ_j+Q_jb_j^L\right)+2r^{-1}\sum_{j=1}^db_j^SQ_j
\end{split}
\end{equation*}
if (\ref{eq:1.7}) holds, and 
\begin{equation*}
\begin{split}
\mathcal{P}(h)&=(\mathcal{D}_r+\sigma)^2+r^{-2}\Lambda-E+h^2(d-1)(d-3)(2r)^{-2}+\widetilde V_L+\widetilde V_S\\
&-\sigma^2-2\sigma \sum_{j=1}^dw_jb_j^S+2\sum_{j=1}^dw_jb_j^S(\mathcal{D}_r+\sigma)\\
&+r^{-1}\sum_{j=1}^d\left(b_j^LQ_j+Q_jb_j^L\right)+2r^{-1}\sum_{j=1}^db_j^SQ_j
\end{split}
\end{equation*}
if (\ref{eq:1.8}) holds. 
Since the function $\varphi$
depends only on the variable $r$, in the first case we get
\begin{equation*}
\begin{split}
\mathcal{P}_\varphi(h)&=\mathcal{D}_r^2+r^{-2}\Lambda-E+h^2(d-1)(d-3)(2r)^{-2}-
2i\varphi'{\cal D}_r+W_L+W_S\\ 
&+2\sum_{j=1}^dw_jb_j\mathcal{D}_r+r^{-1}\sum_{j=1}^d\left(b_j^LQ_j
+Q_jb_j^L\right)+2r^{-1}\sum_{j=1}^d b_j^SQ_j
\end{split}
\end{equation*}
where
\begin{align*}
&W_L=\widetilde V_L-\varphi'^2,\\
&W_S=\widetilde V_S+h\varphi''-2i\varphi'\sum_{j=1}^dw_jb_j+ih\sum_{j=1}^dw_j\partial_r(b_j^L).
\end{align*}
In the second case we get
\begin{equation*}
\begin{split}
\mathcal{P}_\varphi(h)&=(\mathcal{D}_r+\sigma)^2+r^{-2}\Lambda-E+h^2(d-1)(d-3)(2r)^{-2}-
2i\varphi'({\cal D}_r+\sigma)+W_L+W_S\\ 
&+2\sum_{j=1}^dw_jb_j^S(\mathcal{D}_r+\sigma)+r^{-1}\sum_{j=1}^d\left(b_j^LQ_j
+Q_jb_j^L\right)+2r^{-1}\sum_{j=1}^d b_j^SQ_j
\end{split}
\end{equation*}
where
\begin{align*}
&W_L=\widetilde V_L-\sigma^2-\varphi'^2,\\
&W_S=\widetilde V_S+h\varphi''-2(i\varphi'+\sigma)\sum_{j=1}^dw_jb_j^S.
\end{align*}
Set $g=e^{\varphi/h}f$. Since $(\zeta\partial_r)^kf\in H^1(\mathbb{R}^d)$, $k=0,1$, we have 
\begin{equation}\label{eq:3.6}
(\zeta\partial_r)^kg\in H^1(\mathbb{R}^d),\quad k=0,1.
\end{equation} 
Set
\begin{equation*}
P_\varphi(h):=e^{\varphi/h}P(h)e^{-\varphi/h}=P(h)+2h\nabla\varphi\cdot\nabla -|\nabla\varphi|^2+h\Delta\varphi-2ih\nabla\varphi\cdot b.
\end{equation*}
Then $v_1\in D(P_\varphi(h))$ if $P_\varphi(h)v_1\in L^2$ satisfies the identity
\begin{equation}\label{eq:3.7}
\left\langle P_\varphi(h)v_1,v_2\right\rangle_{L^2}=\mathcal{Q}(v_1,v_2),\quad\forall v_2\in H^1,
\end{equation}
where the quadratic form $\mathcal{Q}(v_1,v_2)$ is defined for $v_1,v_2\in H^1$ by
\begin{equation}\label{eq:3.8}
\begin{split}
\mathcal{Q}(v_1,v_2)&=\left\langle (ih\nabla+b)v_1,(ih\nabla+b)v_2\right\rangle_{L^2}\\
&+\left\langle (2h\nabla\varphi\cdot\nabla -|\nabla\varphi|^2+h\Delta\varphi-2ih\nabla\varphi\cdot b+V)v_1,v_2\right\rangle_{L^2}.
\end{split}
\end{equation}
Indeed, since
\begin{equation*}
\Delta\varphi=\varphi''(r)+\frac{d-1}{r}\varphi'(r),\quad r\neq 0,\,1, \, a,
\end{equation*}
by Poincar\'e's inequality $\|r^{-1}v\|_{L^2}\lesssim \|\nabla v\|_{L^2}$, we have
\begin{equation*}
\|\Delta\varphi v_1\|_{L^2}\lesssim \|v_1\|_{H^1}.
\end{equation*}
 Therefore, $\mathcal{Q}(v_1,v_2)$ is well defined for all $v_1,v_2\in H^1$.
 
We will derive the estimate (\ref{eq:3.1}) from the following

\begin{prop} \label{3.3} 
For all functions $g$ satisfying (\ref{eq:3.6}) we have the estimate
\begin{equation}\label{eq:3.9}
\begin{split}
&C\left\|\langle x\rangle^{-s}g\right\|_{L^2}^2+C\left\|\langle x\rangle^{-s}(\mathcal{D}_r^\sharp+\sigma)g\right\|_{L^2}^2\\
&\le -2h^{-1}{\rm Im}\left(\mathcal{Q}(g,\mu
(\mathcal{D}_r^\sharp+\sigma)g)+\left\langle(-E\pm i\varepsilon)g,\mu
(\mathcal{D}_r^\sharp+\sigma)g\right\rangle_{L^2}\right)\\
&+\mathcal{O}(\varepsilon)h^{-1}a^2\left(\tau\|g\|_{L^2}^2+\tau^{-1}\|\mathcal{D}_rg\|_{L^2}^2\right)
\end{split}
\end{equation}
with some constant $C>0$ independent of $\varepsilon$, $h$, $\tau$, $a$ and $g$, where $\mathcal{D}_r^\sharp=\mathcal{D}_r+ih(d-1)(2r)^{-1}$.
\end{prop}

{\it Proof.} Note that, since $\mu(r)=\mathcal{O}(\zeta(r))$
for $0<r\ll 1$, in view of
(\ref{eq:3.8}), the quadratic form in the right-hand side of (\ref{eq:3.9}) is well defined for all functions $g$ satisfying (\ref{eq:3.6}).

We will first prove (\ref{eq:3.9}) for $g\in H^m(\mathbb{R}^d)$ for every integer $m\ge 1$. 
In particular, we have
$g\in C^\infty(\mathbb{R}^d)$. 
Then 
$u(r,w):=r^{(d-1)/2}g(rw)\in C^\infty(\mathbb{R}^+\times\mathbb{S}^{d-1})$. This fact guarantees that the scalar products below are well defined. 
Set
\begin{equation*}
\mathcal{P}_\varphi^\flat(h):=\mathcal{P}_\varphi(h)-r^{-1}\sum_{j=1}^dQ_jb_j^L.
\end{equation*}
Clearly, we have
\begin{equation*}
\mathcal{P}_\varphi^\flat(h)u(r,\cdot)\in L^2(\mathbb{S}^{d-1})
\end{equation*}
for all $r>0$, $r\neq 1$, $r\neq a$. Then, given a function $v\in H^1(\mathbb{S}^{d-1})$, we set
\begin{equation*}
\mathcal{L}(u,v):=\langle \mathcal{P}_\varphi^\flat(h)u,v\rangle+r^{-1}\sum_{j=1}^d
\langle b_j^Lu,Q_jv\rangle.
\end{equation*}
We now define the function
\begin{equation*}
\begin{split}
F(r)=&-\left\langle (r^{-2}\Lambda+h^2(d-1)(d-3)(2r)^{-2}-E+W_L)u(r,\cdot),u(r,\cdot)\right\rangle\\
&+\|(\mathcal{D}_r+\sigma)u(r,\cdot)\|^2-2r^{-1}\sum_{j=1}^d{\rm Re}\left\langle b_j^LQ_ju(r,\cdot),u(r,\cdot)\right\rangle
\end{split}
\end{equation*}
for $r>0$, $r\neq 1$, $r\neq a$. We are going to compute the first derivative of $F$. To this end, 
observe that in the first case considered above, we have
\begin{equation*}
\begin{split}
\frac{h}{2}\frac{d}{dr}\|\mathcal{D}_ru\|^2&={\rm Im}\left\langle\mathcal{D}_r^2u,\mathcal{D}_ru\right\rangle
={\rm Im}\,\mathcal{L}(u,\mathcal{D}_ru)\\
&-{\rm Im}\left\langle (r^{-2}\Lambda+h^2(d-1)(d-3)(2r)^{-2}-E+W_L+W_S)u,\mathcal{D}_ru\right\rangle\\
&+2\varphi'\|\mathcal{D}_ru\|^2-2\sum_{j=1}^d{\rm Im}\left\langle w_jb_j\mathcal{D}_ru,\mathcal{D}_ru\right\rangle\\
&-\sum_{j=1}^d{\rm Im}\left\langle r^{-1}(b_j^LQ_j+Q_jb_j^L)u,\mathcal{D}_ru\right\rangle-2\sum_{j=1}^d{\rm Im}\left\langle r^{-1} b_j^SQ_ju,\mathcal{D}_ru\right\rangle,
\end{split}
\end{equation*}
while in the second case we have
\begin{equation*}
\begin{split}
\frac{h}{2}\frac{d}{dr}\|(\mathcal{D}_r+\sigma)u\|^2&={\rm Im}\left\langle(\mathcal{D}_r+\sigma)^2u,
(\mathcal{D}_r+\sigma)u\right\rangle
={\rm Im}\,\mathcal{L}(u,(\mathcal{D}_r+\sigma)u)\\
&-{\rm Im}\left\langle (r^{-2}\Lambda+h^2(d-1)(d-3)(2r)^{-2}-E+W_L+W_S)u,(\mathcal{D}_r+\sigma)u\right\rangle\\
&+2\varphi'\|(\mathcal{D}_r+\sigma)u\|^2-2\sum_{j=1}^d{\rm Im}\left\langle w_jb_j^S
(\mathcal{D}_r+\sigma)u,(\mathcal{D}_r+\sigma)u\right\rangle\\
&-\sum_{j=1}^d{\rm Im}\left\langle r^{-1}(b_j^LQ_j+Q_jb_j^L)u,(\mathcal{D}_r+\sigma)u\right\rangle\\
&-2\sum_{j=1}^d{\rm Im}\left\langle r^{-1} b_j^SQ_ju,(\mathcal{D}_r+\sigma)u\right\rangle.
\end{split}
\end{equation*}
On the other hand, we have the identities
\begin{equation*}
\begin{split}
&{\rm Im}\left\langle w_jb_j\mathcal{D}_ru,\mathcal{D}_ru\right\rangle=0,\\
&{\rm Re}\frac{d}{dr}\left\langle b_j^LQ_ju,u\right\rangle={\rm Re}\left\langle 
\partial_r(b_j^L)Q_ju,u\right\rangle\\
&-h^{-1}{\rm Im}\left\langle(b_j^LQ_j+Q_jb_j^L)u,\mathcal{D}_ru\right\rangle
\end{split}
\end{equation*}
in the first case, and
\begin{equation*}
\begin{split}
&{\rm Im}\left\langle w_jb_j^S(\mathcal{D}_r+\sigma)u,(\mathcal{D}_r+\sigma)u\right\rangle=0,\\
&{\rm Re}\frac{d}{dr}\left\langle b_j^LQ_ju,u\right\rangle={\rm Re}\left\langle 
\partial_r(b_j^L)Q_ju,u\right\rangle\\
&-h^{-1}{\rm Im}\left\langle(b_j^LQ_j+Q_jb_j^L)u,(\mathcal{D}_r+\sigma)u\right\rangle-
{\rm Re}\left\langle b_j^Lu,[q_j,\sigma]u\right\rangle
\end{split}
\end{equation*}
in the second case. 
Using the above identities we find that in both cases the first derivative of $F$ is given by
\begin{equation*}
\begin{split}
F'(r)&=\frac{2}{r}\left\langle r^{-2}(\Lambda+(h/2)^2(d-1)(d-3))u,u\right\rangle
-W'_L\|u\|^2\\
&+2r^{-2}\sum_{j=1}^d{\rm Re}\left\langle b_j^LQ_ju,u\right\rangle
-2r^{-1}\sum_{j=1}^d{\rm Re}\left\langle\partial_r(b_j^L)Q_ju,u\right\rangle
+2r^{-1}\sum_{j=1}^d{\rm Re}\left\langle b_j^Lu,[q_j,\sigma]u\right\rangle\\
&+2h^{-1}{\rm Im}\,\mathcal{L}(u,(\mathcal{D}_r+\sigma)u)+4h^{-1}
\varphi'\|(\mathcal{D}_r+\sigma)u\|^2-2h^{-1}{\rm Im}\,\left\langle\widetilde W_Su,({\cal D}_r+\sigma)u\right\rangle,
\end{split}
\end{equation*}
where
\begin{equation*}
\widetilde W_S=W_S+2r^{-1}\sum_{j=1}^d b_j^SQ_j.
\end{equation*}
Thus we obtain the identity
\begin{equation*}
\begin{split}
(\mu F)'(r)&=\mu'(r)F(r)+\mu(r)F'(r)\\
&=(2r^{-1}\mu-\mu')\left\langle r^{-2}(\Lambda+(h/2)^2(d-1)(d-3))u,u\right\rangle\\ 
&+(E\mu'-(\mu W_L)')\|u\|^2+(\mu'+4h^{-1}\varphi'\mu)\|(\mathcal{D}_r+\sigma)u\|^2\\
&-2(r^{-1}\mu)'\sum_{j=1}^d{\rm Re}\left\langle b_j^LQ_ju,u\right\rangle-2r^{-1}\mu
\sum_{j=1}^d{\rm Re}\left\langle\partial_r(b_j^L)Q_ju,u\right\rangle\\
&+2r^{-1}\mu \sum_{j=1}^d{\rm Re}\left\langle b_j^Lu,[q_j,\sigma]u\right\rangle-2h^{-1}\mu{\rm Im}\,\left\langle\widetilde W_Su,
(\mathcal{D}_r+\sigma)u\right\rangle\\
&+2h^{-1}\Psi_{\pm}(r)\mp 2\varepsilon h^{-1}\mu{\rm Re}\,\left\langle u,(\mathcal{D}_r+\sigma)u\right\rangle,
\end{split}
\end{equation*}
where 
\begin{equation*}
\Psi_{\pm}(r)={\rm Im}\,\mathcal{L}(u(r,\cdot),\mu(r)
(\mathcal{D}_r+\sigma)u(r,\cdot))+{\rm Im}\left\langle\pm i\varepsilon u(r,\cdot),\mu(r)
(\mathcal{D}_r+\sigma)u(r,\cdot)\right\rangle.
\end{equation*}
In view of (\ref{eq:3.5}), we have the estimate
\begin{equation}\label{eq:3.10}
\left\|\widetilde W_Su\right\|^2\lesssim B_1\left\|u\right\|^2+B_2\left\|r^{-1}\Lambda^{1/2}u\right\|^2,
\end{equation}
where 
\begin{align*}
&B_1(r)=(r+1)^{-2-\rho}+(r+1)^{-2-\rho}\varphi'^2+h^2\varphi''^2+h^2r^{-2}(r+1)^{-2-\rho},\\
&B_2(r)=(r+1)^{-2-\rho}.
\end{align*}
We also have
\begin{equation}\label{eq:3.11}
\begin{split}
\left|\left\langle b_j^LQ_ju,u\right\rangle\right|&\lesssim B_3\|u\|\left\|\Lambda^{1/2}u\right\|+hB_3\|u\|^2\\
 &\lesssim \left(\gamma_1^{-1}(rB_3)^2+hB_3\right)\|u\|^2+\gamma_1\left\|r^{-1}\Lambda^{1/2}u\right\|^2,
 \end{split}
\end{equation}
\begin{equation}\label{eq:3.12}
\left|\left\langle\partial_r(b_j^L)Q_ju,u\right\rangle\right|\lesssim \left(\gamma_2^{-1}(rB_4)^2+hB_4\right)\|u\|^2+\gamma_2\left\|r^{-1}\Lambda^{1/2}u\right\|^2,
\end{equation}
where $\gamma_1$ and $\gamma_2$ are arbitrary positive functions of $r$ to be fixed below and 
\begin{equation*}
B_3(r)=(r+1)^{-\rho},\quad B_4(r)=(r+1)^{-1-\rho}.
\end{equation*}
Write now (\ref{eq:3.3}) in the form
\begin{equation*}
w_j\partial_r+r^{-1}q_j(w,\partial_w)=\partial_{x_j}-\frac{d-1}{4r^2}\frac{\partial r^2}{\partial x_j} 
=\partial_{x_j}-(d-1)(2r)^{-1}w_j.
\end{equation*}
Thus we obtain the identity
\begin{equation*}
[q_j,\sigma]=[r^{-1}q_j,r\sigma]=\partial_{x_j}(x\cdot b^L(x))-w_j\partial_r(r\sigma).
\end{equation*}
Hence, in view of (\ref{eq:1.4}) and (\ref{eq:1.8}), we get the bound
\begin{equation*}
\left|[q_j,\sigma]\right|\lesssim 1,
\end{equation*}
which implies 
\begin{equation}\label{eq:3.13}
\left|\left\langle b_j^Lu,[q_j,\sigma]u\right\rangle\right|\lesssim B_5\|u\|^2,
\end{equation}
where $B_5(r)=(r+1)^{-\rho}$. 
Since $2r^{-1}\mu-\mu'>0$ and $d\ge 3$, by (\ref{eq:3.10}), (\ref{eq:3.11}), (\ref{eq:3.12}) and (\ref{eq:3.13}), we get the inequality
\begin{equation*}
\begin{split}
(\mu F)'(r)&\ge (2r^{-1}\mu-\mu')\left\langle r^{-2}\Lambda u,u\right\rangle\\
&+(E\mu'-(\mu W_L)')\|u\|^2+(\mu'+4h^{-1}\varphi'\mu)\|(\mathcal{D}_r+\sigma)u\|^2\\
&-2|(r^{-1}\mu)'|\sum_{j=1}^d\left|\left\langle b_j^LQ_ju,u\right\rangle\right|-2r^{-1}\mu\sum_{j=1}^d
\left|\left\langle\partial_r(b_j^L)Q_ju,u\right\rangle\right|
-2r^{-1}\mu\sum_{j=1}^d\left|\left\langle b_j^Lu,[q_j,\sigma]u\right\rangle\right|\\ 
&-2h^{-2}\mu^2(\mu'+4h^{-1}\varphi'\mu)^{-1}\left\|\widetilde W_Su\right\|^2
-\frac{1}{2}(\mu'+4h^{-1}\varphi'\mu)\|(\mathcal{D}_r+\sigma)u\|^2\\
&+2h^{-1}\Psi_{\pm}(r)-\varepsilon h^{-1}\mu\|u\|\|(\mathcal{D}_r+\sigma)u\|\\
&\ge A_1\|u\|^2+A_2\|r^{-1}\Lambda^{1/2}u\|^2+\frac{\mu'}{2}\|(\mathcal{D}_r+\sigma)u\|^2+
2h^{-1}\Psi_{\pm}(r)-\varepsilon h^{-1}\mu\|u\|\|(\mathcal{D}_r+\sigma)u\|,
\end{split}
\end{equation*}
where 
\begin{equation*}
\begin{split}
A_1(r)&=(\varphi'^2\mu)'+(E-\widetilde V_L)\mu'-C(r+1)^{-1-\rho}\mu-Ch^{-2}B_1\mu^2(\mu'+h^{-1}\varphi'\mu)^{-1}\\
&-C|(r^{-1}\mu)'|\left(\gamma_1^{-1}(rB_3)^2+hB_3\right)-Cr^{-1}\mu\left(\gamma_2^{-1}(rB_4)^2+hB_4+B_5\right),\\
A_2(r)&=2r^{-1}\mu-\mu'-Ch^{-2}B_2\mu^2(\mu'+h^{-1}\varphi'\mu)^{-1}\\
&-C\gamma_1|(r^{-1}\mu)'|-C\gamma_2r^{-1}\mu.
\end{split}
\end{equation*}
 We now choose the functions $\gamma_1$ and $\gamma_2$ so that
\begin{equation*}
C\gamma_1|(r^{-1}\mu)'|=C\gamma_2r^{-1}\mu=\frac{1}{4}(2r^{-1}\mu-\mu').
\end{equation*}
Thus we conclude that the above inequality holds with
\begin{equation*}
\begin{split}
A_1(r)&=(\varphi'^2\mu)'+(E-\widetilde V_L)\mu'-Cr^{-1}(r+1)^{-\rho}\mu-Ch^{-2}B_1\mu^2(\mu'+h^{-1}\varphi'\mu)^{-1}\\
&-C(rB_3)^2M_1-ChB_3M_2,\\
2A_2(r)&=2r^{-1}\mu-\mu'-Ch^{-2}B_2\mu^2(\mu'+h^{-1}\varphi'\mu)^{-1},
\end{split}
\end{equation*}
with a new constant $C>0$, where
\begin{equation*}
\begin{split}
M_1(r)&=(2r^{-1}\mu-\mu')^{-1}\left(|(r^{-1}\mu)'|^2+(r(r+1))^{-2}\mu^2\right),\\
M_2(r)&=|(r^{-1}\mu)'|+(r(r+1))^{-1}\mu.
\end{split}
\end{equation*} 
 Now we will use Lemmas \ref{2.1} and \ref{2.3} to bound the functions $A_1$ and $A_2$ from below.
 Observe first that $M_1(r)=M_2(r)=1$ for $r<1$. For $1<r<a$ we have the bounds
 \begin{equation*}
 M_1(r)\lesssim (r+1)^{-1+\ell},\quad M_2(r)\lesssim 1,
 \end{equation*}
 while for $r>a$ we have the bounds
 \begin{equation*}
 M_1(r)\lesssim (r+1)^{-3+2s}\mu'(r),\quad M_2(r)\lesssim (r+1)^{-2+2s}\mu'(r).
 \end{equation*}
 From the properties of the functions $\varphi$ and $\mu$ we also deduce
 \begin{equation*}
 B_1(r)\mu^2(\mu'+h^{-1}\varphi'\mu)^{-1}\lesssim
 \tau^2h\varphi'(r)^{-1}\mu(r)+h^2r^{-2}\mu(r)^2\mu'(r)^{-1}\lesssim \tau^2h
 \end{equation*}
 for $0<r<1$,
 \begin{equation*}
 \begin{split}
  B_1(r)\mu^2(\mu'+h^{-1}\varphi'\mu)^{-1}&\lesssim
  h B_1(r)\mu(r)\varphi'(r)^{-1}\\
  &\lesssim
  h\tau^{-1}(r+1)^{2}B_1(r)\mu'(r)\\
  &\lesssim  h\tau^{-1}(r+1)^{-\rho}\mu'(r)+h\tau(r+1)^{-1-\rho}
  \end{split}
 \end{equation*}
 for $1<r<a$,
 \begin{equation*}
 \begin{split}
 B_1(r)\mu^2(\mu'+h^{-1}\varphi'\mu)^{-1}&\lesssim
  h B_1(r)\mu(r)\varphi'(r)^{-1}\\
  &\lesssim
  h (r+1)^{2s}B_1(r)\varphi'(r)^{-1}\mu'(r)\\
  &\lesssim
   h\tau^{-1}a^{-\kappa}(r+1)^{2s+\kappa-1-\rho}\mu'(r)\\
   &+h\tau a^{\kappa}(r+1)^{2s-3-\kappa-\rho}\mu'(r)\\
    &+h^3\tau a^{\kappa}\tau(r+1)^{2s-3-\kappa}\mu'(r)\\
    &\lesssim
   h\tau^{-1}(r+1)^{2s-1-\rho}\mu'(r)+h\tau(r+1)^{2s-3}\mu'(r)
   \end{split}
  \end{equation*}
  for $r>a$ if we take $1/2<s<(1+\rho)/2$ and $0<\kappa<1+\rho-2s$. 
  
  Thus, for $0<r<1$ we obtain
\begin{equation*}
\begin{split}
A_1(r)&\ge E\mu'(r)+C^\flat_\nu\tau^2-Ch^{-2}-Ch^{-1}\tau,\quad C^\flat_\nu>0,\\
A_2(r)&\ge 2^{-1}-Ch^{-1}\tau^{-1}.
\end{split}
\end{equation*}
We take $\tau=\tau_0h^{-1}$, where $\tau_0\gg 1$ is independent of $h$.
It is clear that if $\tau_0$ is big enough we can arrange that
 \begin{equation*}
A_1(r)\ge E\mu'(r),\quad A_2(r)\ge 0, \quad\mbox{for}\quad 0<r<1.
 \end{equation*}
For $1<r<a$ we have $r^{-1}\mu(r)<\mu'(r)$ and 
\begin{equation}\label{eq:3.14}
(rB_3)^2M_1(r)\lesssim (r+1)^{\ell-2\rho}\mu'(r)\lesssim (r+1)^{-\rho}\mu'(r)
\end{equation}
if $\ell<\rho$. Hence
\begin{equation*}
\begin{split}
A_1(r)&\ge \tau^2 C_\nu(r+1)^{-1-\nu}+(E-\widetilde V_L)\mu'-C(r+1)^{-\rho}\mu'\\
&-Ch^{-1}\tau^{-1}(r+1)^{-\rho}\mu'-Ch\tau(r+1)^{-2}\mu'-Ch^{-1}\tau(r+1)^{-1-\rho},\\
A_2(r)&\ge 2^{-1}\ell r^{1-\ell}-Ch^{-1}\tau^{-1}r^{1-\rho}.
\end{split}
\end{equation*}
 Since $\widetilde V_L$ satisfies (\ref{eq:1.1}),
there exists $\lambda_0\gg 1$, independent of $h$ and $\tau_0$, such that
\begin{equation}\label{eq:3.15}
\widetilde V_L(rw)+C(r+1)^{-\rho}\le E/4\quad\mbox{for}\quad r\ge\lambda_0.
\end{equation}
On the other hand, for $r<\lambda_0$, we have
\begin{equation}\label{eq:3.16}
\widetilde V_L\mu'(r)+C(r+1)^{-\rho}\mu'(r)\lesssim (r+1)^{-1-\nu}.
\end{equation}
Similarly, there exists $\lambda\gg 1$, independent of $h$ but depending on $\tau_0$, such that
\begin{equation*}
C\tau_0(r+1)^{-2}\le E/4\quad\mbox{for}\quad r\ge\lambda,
\end{equation*}
and
\begin{equation*}
C\tau_0(r+1)^{-2}\mu'(r)\lesssim (r+1)^{-1-\nu}\quad\mbox{for}\quad r<\lambda.
\end{equation*}
Therefore, if we choose $\nu<\rho$, $\ell<\rho$, taking $\tau_0$ big enough and $h$ small enough we can arrange that
\begin{equation*}
A_1(r)\ge 3^{-1}E\mu'(r),\quad A_2(r)\ge 0, \quad\mbox{for}\quad 1<r<a.
\end{equation*}
In view of (\ref{eq:1.1}) with $\widetilde V_L$, for $r>a$ we have
\begin{equation*}
\begin{split}
A_1(r)&\ge 2^{-1}E\mu'(r)-C\mu'(r)\left(\tau^2a^{-3+2s}+a^{-\rho+2s-1}\right),\\
A_2(r)&\ge 2^{-1}a^2(r+1)^{-1}\left(1-Ch^{-1}\tau^{-1}a^{-1}\right).
\end{split}
\end{equation*}
Taking $a=h^{-m}$ with $m>0$ big enough and $h$ small enough we can arrange that
\begin{equation*}
A_1(r)\ge 3^{-1}E\mu'(r),\quad A_2(r)\ge 0, \quad\mbox{for}\quad r>a.
\end{equation*}
By the above inequalities we obtain
\begin{equation}\label{eq:3.17}
(\mu F)'(r)\ge \frac{E}{3}\mu'\|u\|^2+\frac{\mu'}{2}\|(\mathcal{D}_r+\sigma)u\|^2+
2h^{-1}\Psi_{\pm}(r)
-\varepsilon h^{-1}\mu\|u\|\|(\mathcal{D}_r+\sigma)u\|.
\end{equation}
Integrating (\ref{eq:3.17}) with respect to $r$ and using that 
\begin{equation*}
\int_0^\infty(\mu F)'dr=0,
\end{equation*}
 we obtain the estimate
\begin{equation}\label{eq:3.18}
\begin{split}
&\frac{E}{3}\int_0^\infty\mu'\|u\|^2dr+\frac{1}{2}\int_0^\infty\mu'\|(\mathcal{D}_r+\sigma)u\|^2dr\\
&\le -2h^{-1}\int_0^\infty\Psi_{\pm}(r)dr
+\varepsilon h^{-1}\int_0^\infty\mu\|u\|\|(\mathcal{D}_r+\sigma)u\|dr.
\end{split}
\end{equation}
In what follows we will use that $\mu=\zeta\widetilde\mu$ with a Lipschitz function $\widetilde\mu$ satisfying
 $|\partial_r^\ell\widetilde \mu|\lesssim a^2$, $\ell=0,1$.
 In view of (\ref{eq:2.3}), we get from (\ref{eq:3.18}),
\begin{equation}\label{eq:3.19}
\begin{split}
&C\int_0^\infty(r+1)^{-2s}\left(\|u\|^2+\|(\mathcal{D}_r+\sigma)u\|^2\right)dr\\
&\le -2h^{-1}\int_0^\infty\Psi_{\pm}(r)dr+\varepsilon h^{-1}a^2\int_0^\infty\left(\tau\|u\|^2+\tau^{-1}\zeta^2\|\mathcal{D}_ru\|^2\right)dr
\end{split}
\end{equation}
with some constant $C>0$ independent of $h$, $\tau$, $a$ and $\varepsilon$. Since
\begin{equation*}
u=r^{(d-1)/2}g,\quad \mathcal{D}_ru=r^{(d-1)/2}\mathcal{D}_r^\sharp g,
\end{equation*}
we have
\begin{equation}\label{eq:3.20}
\zeta^2\|\mathcal{D}_ru\|^2\lesssim \left\|r^{(d-1)/2}\mathcal{D}_rg\right\|^2+\left\|r^{(d-1)/2}g\right\|^2.
\end{equation} 
Let us see now that we have the identity
\begin{equation}\label{eq:3.21}
\int_0^\infty\Psi_{\pm}(r)dr={\rm Im}\,\mathcal{Q}(g,\mu
(\mathcal{D}_r^\sharp+\sigma)g)+{\rm Im}\left\langle(-E\pm i\varepsilon)g,\mu
(\mathcal{D}_r^\sharp+\sigma)g\right\rangle_{L^2}.
\end{equation}
Set $g_1:=\mu
(\mathcal{D}_r^\sharp+\sigma)g$ and 
\begin{equation*}
P^\flat(h):=P(h)-r^{-1}\sum_{j=1}^dQ_jb_j^L,
\end{equation*}
\begin{equation*}
P_\varphi^\flat(h):=P_\varphi(h)-r^{-1}\sum_{j=1}^dQ_jb_j^L=P^\flat(h)+2h\nabla\varphi\cdot\nabla -|\nabla\varphi|^2+h\Delta\varphi-2ih\nabla\varphi\cdot b.
\end{equation*}
We have
\begin{equation}\label{eq:3.22}
\begin{split}
\int_0^\infty\mathcal{L}(r^{(d-1)/2}g,r^{(d-1)/2}g_1)dr&=\left\langle(P_\varphi^\flat(h)-E)g,g_1\right\rangle_{L^2}
+\sum_{j=1}^d\left\langle r^{-1}b_j^Lg,Q_jg_1\right\rangle_{L^2}\\
&=\left\langle(P^\flat(h)-E)g,g_1\right\rangle_{L^2}
+\sum_{j=1}^d\left\langle r^{-1}b_j^Lg,Q_jg_1\right\rangle_{L^2}\\
&+\left\langle(2h\nabla\varphi\cdot\nabla -|\nabla\varphi|^2+h\Delta\varphi-2ih\nabla\varphi\cdot b)g,g_1\right\rangle_{L^2}\\
&=\mathcal{Q}(g,g_1)-E\langle g,g_1\rangle_{L^2}.
\end{split}
\end{equation}
It is easy to see that (\ref{eq:3.21}) follows from (\ref{eq:3.22}). 

By (\ref{eq:3.19}), (\ref{eq:3.20}) and (\ref{eq:3.21}), we obtain that the inequality (\ref{eq:3.9})
holds for all $g\in H^m$, $\forall m\in\mathbb{N}$. By a density argument, we will show that (\ref{eq:3.9})
holds for all $g$ satisfying (\ref{eq:3.6}). Let $\Phi\in C_0^\infty(\mathbb{R}^d)$, $\Phi\ge 0$, be such that
$\int \Phi(x)dx=1$. Set $\Phi_\epsilon(x)=\epsilon^{-d}\Phi(x/\epsilon)$, $0<\epsilon\ll 1$. Then 
$g_\epsilon:=\Phi_\epsilon*g\in H^1$ uniformly in $\epsilon$. Moreover, $\|g-g_\epsilon\|_{H^1}\to 0$ as $\epsilon\to 0$. 
It also follows from Young's inequality
that $g_\epsilon\in H^m$ for all integers $m\ge 1$ and $\|g_\epsilon\|_{H^m}=\mathcal{O}(\epsilon^{1-m})$. 
To approximate the function $\zeta\partial_rg$ we will use the Friedrichs lemma (see Lemma 17.1.5 of \cite{kn:H}):

\begin{lemma} \label{3.4} 
Let $g\in L^2(\mathbb{R}^d)$ and let $\omega$ be a Lipschitz function satisfying the inequality 
\begin{equation}\label{eq:3.23}
 |\omega(x)-\omega(y)|\le C|x-y|,\quad\forall x,y\in \mathbb{R}^d.
\end{equation}
Then
\begin{equation}\label{eq:3.24}
 \|\omega(\Phi_\epsilon*\partial_{x_j}g)-\Phi_\epsilon*(\omega\partial_{x_j} g)\|_{L^2}\to 0\quad\mbox{as}\quad \epsilon\to 0.
\end{equation}
\end{lemma}

We have 
\begin{equation*}
\zeta\partial_r=\sum_{\nu=1}^d\omega_\nu(x)\partial_{x_\nu},\quad \omega_\nu=x_\nu(1+|x|^2)^{-1/2}.
\end{equation*}
It is easy to see that the functions $\omega_\nu$ satisfy (\ref{eq:3.23}). Therefore, 
by (\ref{eq:3.24}) we get  
\begin{equation*}
\|\zeta\partial_r(\Phi_\epsilon*g)-\Phi_\epsilon*(\zeta\partial_rg)\|_{L^2}\to 0\quad\mbox{as}\quad\epsilon\to 0,
\end{equation*}
which implies 
$\|\zeta\partial_rg-\zeta\partial_rg_\epsilon\|_{L^2}\to 0$ as $\epsilon\to 0$. 
Clearly, this still holds with $g$ replaced by $\partial_{x_j}g$. This in turn imples that the same limit holds also
in the $H^1$ norm.  Indeed, using that
\begin{equation*}
[\zeta\partial_r,\partial_{x_j}]=-\sum_{\nu=1}^d\frac{\partial \omega_\nu}{\partial x_j}\partial_{x_\nu},
\end{equation*}
we get
\begin{equation*}
\|\partial_{x_j}\zeta\partial_r(g-g_\epsilon)\|_{L^2}\lesssim 
\sum_{k=0}^1\sum_{\nu=1}^d\|(\zeta\partial_r)^k\partial_{x_\nu}(g-g_\epsilon)\|_{L^2},
\end{equation*}
which implies $\|\zeta\partial_rg-\zeta\partial_rg_\epsilon\|_{H^1}\to 0$ as $\epsilon\to 0$, as desired.

 In view of (\ref{eq:3.8}), taking into account that $\mu=\zeta\widetilde\mu$, one can easily get the estimate
\begin{equation}\label{eq:3.25}
\begin{split}
&\left|\mathcal{Q}(g,\mu
(\mathcal{D}_r^\sharp+\sigma)g)+\left\langle(-E\pm i\varepsilon)g,\mu
(\mathcal{D}_r^\sharp+\sigma)g\right\rangle_{L^2}\right.\\
&\left.-\mathcal{Q}(g_\epsilon,\mu
(\mathcal{D}_r^\sharp+\sigma)g_\epsilon)-\left\langle(-E\pm i\varepsilon)g_\epsilon,\mu
(\mathcal{D}_r^\sharp+\sigma)g_\epsilon\right\rangle_{L^2}\right|\\
&\le C\left(\|g\|_{H^1}+\|\zeta\partial_rg\|_{H^1}\right)\left(\|g-g_\epsilon\|_{H^1}
+\|\zeta\partial_rg-\zeta\partial_rg_\epsilon\|_{H^1}\right)
\end{split}
\end{equation}
with a constant $C>0$ independent of $g$ and $\epsilon$. 
By Poincar\'e's inequality $\|r^{-1}v\|_{L^2}\lesssim \|\nabla v\|_{L^2}$, we also have 
\begin{equation}\label{eq:3.26}
\left\|(\mathcal{D}_r^\sharp+\sigma)(g-g_\epsilon)\right\|_{L^2}\le 
C\|g-g_\epsilon\|_{H^1}.
\end{equation}
Hence the left-hand sides of (\ref{eq:3.25}) and (\ref{eq:3.26}) tend to zero as $\epsilon\to 0$.
Thus we conclude that if the estimate (\ref{eq:3.9}) holds with $g_\epsilon$ uniformly in $\epsilon$, then it also holds with $g$, as desired.
\eproof

In addition to (\ref{eq:3.6}) suppose that $g\in D(P_\varphi(h))$. 
 We are going to use the identity (\ref{eq:3.7}) with
\begin{equation*}
v_1=g,\quad v_2=\mu(\mathcal{D}_r^\sharp+\sigma)g.
\end{equation*}
Then the first term in the right-hand side of (\ref{eq:3.9}) can be written in the form
\begin{equation*}
-2h^{-1}{\rm Im}\left\langle(P_\varphi(h)-E\pm i\varepsilon)g,\mu
(\mathcal{D}_r^\sharp+\sigma)g\right\rangle_{L^2},
\end{equation*}
so it is bounded from above by
\begin{equation*}
\mathcal{O}(\gamma^{-1})a^2h^{-2}\left\|\langle x\rangle^{s}(P_\varphi(h)-E\pm i\varepsilon)g\right\|_{L^2}^2+\gamma\left\|\langle x\rangle^{-s}(\mathcal{D}_r^\sharp+\sigma)g\right\|_{L^2}^2
\end{equation*}
for any $0<\gamma\le 1$ independent of $h$ and $a$. Therefore, taking $\gamma$ as small as needed, we deduce from (\ref{eq:3.9}),
\begin{equation}\label{eq:3.27}
 \begin{split}
\|\langle x\rangle^{-s}g\|_{L^2}&\le Cah^{-1}\|\langle x\rangle^{s}(P_\varphi(h)-E\pm i\varepsilon)g\|_{L^2}\\ 
&+Ca\left(\tau h^{-1}\varepsilon\right)^{1/2}\|g\|_{L^2}\\
&+Ca\left(\tau^{-1}h^{-1}\varepsilon \right)^{1/2}\|{\cal D}_rg\|_{L^2}
\end{split}
\end{equation}
with some constant $C>0$ independent of $h$, $\tau$, $a$ and $\varepsilon$. 
Since
\begin{equation*}
{\cal D}_r(e^{\varphi/h}f)=e^{\varphi/h}({\cal D}_rf+i\tau\psi'f),
\end{equation*}
it is easy to see that the estimate (\ref{eq:3.1}) follows from (\ref{eq:3.27}).
\eproof

\section{Carleman estimates when $V_S\equiv 0$ and $b^S\equiv 0$}

In this section we prove the following

\begin{Theorem} \label{4.1}
Under the conditions of Theorem \ref{1.1} with (\ref{eq:1.8}), if $V_S\equiv 0$ and $b^S\equiv 0$, 
 the estimate (\ref{eq:3.1}) holds true with new parameters $a,\tau\gg 1$ independent of $h$.
\end{Theorem}

{\it Proof.} We will make use of the Carleman estimates obtained in the previous section under the condition 
(\ref{eq:1.8}), which in our case take a much simpler form. Indeed, $V_S\equiv 0$, $b^S\equiv 0$ imply that
$\widetilde V_S\equiv 0$ and $\widetilde W_S=h\varphi''$. Hence $B_1=h^2\varphi''^2$ and $B_2=0$. The functions 
$B_3$, $B_4$ and $B_5$ remain unchanged. Therefore the functions $A_1$ and $A_2$ take the form
\begin{equation*}
\begin{split}
A_1(r)&=(\varphi'^2\mu)'+(E-\widetilde V_L)\mu'-Cr^{-1}(r+1)^{-\rho}\mu\\
&-Ch\tau\mu(\psi'')^2(\psi')^{-1}-C(rB_3)^2M_1-ChB_3M_2,\\
2A_2(r)&=2r^{-1}\mu-\mu'>0.
\end{split}
\end{equation*}
For $0<r<a$, $r\neq 1$, in view of (\ref{eq:3.14}), (\ref{eq:3.15}) and (\ref{eq:3.16}), we have
\begin{equation*}
A_1(r)\ge \tau^2C_\nu(r+1)^{-1-\nu}+(E-\widetilde V_L)\mu'-C(r+1)^{-\rho}\mu'-Ch\tau(r+1)^{-2}\mu'\ge E\mu'/3,
\end{equation*}
if $\tau$ is taken big enough, independent of $h$, and $h$ is taken small enough. For $r>a$ we have
\begin{equation*}
A_1(r)\ge 2^{-1}E\mu'-C\mu'\left(\tau^2a^{-3+2s}+a^{-\rho+2s-1}\right)\ge E\mu'/3,
\end{equation*}
if $a$ is taken big enough, independent of $h$. Thus we conclude that the estimate (\ref{eq:3.19}) still holds 
in this case with parameters $a,\tau\gg 1$ independent of $h$. Clearly, this implies the theorem.
\eproof

\section{Carleman estimates for H\"older potentials}

In this section we prove the following

\begin{Theorem} \label{5.1}
Under the conditions of Theorem \ref{1.2}, the estimate (\ref{eq:3.1}) holds true with new parameters $a=h^{-m}$
and $\tau=\tau_0h^{-n+1}$ with some positive constants $m$ and $\tau_0$.
\end{Theorem}

{\it Proof.} 
Let $\phi_0\in C_0^\infty([0,1])$, $\phi_0\ge 0$, be a real-valued function independent of $h$ such that $\int_0^\infty\phi_0(t)dt=1$.
If $V_S$ satisfies (\ref{eq:1.21}), we approximate it by the function
\begin{equation*}
V_{\theta_1}(r,w)=\theta_1^{-1}\int_0^\infty\phi_0((r'-r)/\theta_1)V_S(r',w)dr'=
\int_0^\infty\phi_0(t)V_S(r+\theta_1t,w)dt
\end{equation*}
where $0<\theta_1<1$ is a parameter depending on $h$ to be fixed later on. We have
\begin{equation}\label{eq:5.1}
\begin{split}
|V_S(r,w)-V_{\theta_1}(r,w)|&\le\int_0^\infty\phi_0(t)|V_S(r+\theta_1t,w)-V_S(r,w)|dt\\
&\lesssim\theta_1^\alpha(r+1)^{-3-\rho}\int_0^\infty t^\alpha\phi_0(t)dt\\
&\lesssim\theta_1^\alpha(r+1)^{-3-\rho}.
\end{split}
\end{equation}
Since $V_S$ satisfies (\ref{eq:1.1}), the above bound implies
\begin{equation}\label{eq:5.2}
V_{\theta_1}(r,w)\le p(r)+{\cal O}((r+1)^{-3-\rho}).
\end{equation}
Clearly, $V_{\theta_1}$ is $C^1$ with respect to the variable $r$ and its first derivative is given by
\begin{equation*}
\begin{split}
\partial_rV_{\theta_1}(r,w)&=-\theta_1^{-2}\int_0^\infty\phi_0'((r'-r)/\theta_1)V_S(r',w)dr'\\
&=-\theta_1^{-1}\int_0^\infty\phi_0'(t)V_S(r+\theta_1t,w)dt\\
&=-\theta_1^{-1}\int_0^\infty\phi_0'(t)(V_S(r+\theta_1t,w)-V_S(r,w))dt
\end{split}
\end{equation*}
where we have used that $\int_0^\infty\phi_0'(t)dt=0$. Hence
\begin{equation}\label{eq:5.3}
|\partial_rV_{\theta_1}(r,w)|\lesssim\theta_1^{-1+\alpha}(r+1)^{-3-\rho}\int_0^\infty t^\alpha|\phi_0'(t)|dt
\lesssim\theta_1^{-1+\alpha}(r+1)^{-3-\rho}.
\end{equation}
Let $\phi\in C_0^\infty(\mathbb{R}^d)$, $\phi\ge 0$, be a real-valued function independent of $h$ 
such that $\int_{\mathbb{R}^d}\phi(x)dx=1$.
If $b^S$ satisfies (\ref{eq:1.22}), we approximate it by the function
\begin{equation*}
b_{\theta_2}(x)=\theta_2^{-d}\int_{\mathbb{R}^d}\phi((x'-x)/\theta_2)b^S(x')dx'=
\int_{\mathbb{R}^d}\phi(y)b^S(x+\theta_2y)dy
\end{equation*}
where $0<\theta_2<1$ is a parameter depending on $h$ to be fixed later on. We have
\begin{equation}\label{eq:5.4}
\begin{split}
|b^S(x)-b_{\theta_2}(x)|&\le\int_{\mathbb{R}^d}\phi(y)|b^S(x+\theta_2y)-b^S(x)|dy\\
&\lesssim\theta_2^{\alpha'}(r+1)^{-3/2-\rho_1}\int_{\mathbb{R}^d} |y|^{\alpha'}\phi(y)dy\\
&\lesssim\theta_2^{\alpha'}(r+1)^{-3/2-\rho_1}.
\end{split}
\end{equation}
By (\ref{eq:1.23}) and (\ref{eq:5.4}),
\begin{equation}\label{eq:5.5}
|b_{\theta_2}(rw)|\lesssim (r+1)^{-\rho_2}+(r+1)^{-3/2-\rho_1}.
\end{equation}
Since $\rho_1+\rho_2>3/2$ by assumption, we deduce from (\ref{eq:1.23}), (\ref{eq:5.4}) and (\ref{eq:5.5}),
\begin{equation}\label{eq:5.6}
\begin{split}
&|(b^S-b_{\theta_2})(rw)||(b^L+b_{\theta_2})(rw)|\\
&\lesssim \theta_2^{\alpha'}(r+1)^{-3/2-\rho_1-\rho_2}+\theta_2^{\alpha'}(r+1)^{-3-2\rho_1}\\
&\lesssim \theta_2^{\alpha'}(r+1)^{-3-\rho}
\end{split}
\end{equation}
with some constant $\rho>0$. Furthermore, if $\beta$ is a multi-index such that $|\beta|=1$, we have
\begin{equation*}
\begin{split}
\partial_x^\beta b_{\theta_2}(x)&=-\theta_2^{-d-1}\int_{\mathbb{R}^d}(\partial_x^\beta\phi)((x'-x)/\theta_2)b^S(x')dx'\\
&=-\theta_2^{-1}\int_{\mathbb{R}^d}(\partial_y^\beta\phi)(y)b^S(x+\theta_2y)dy\\
&=-\theta_2^{-1}\int_{\mathbb{R}^d}(\partial_y^\beta\phi)(y)(b^S(x+\theta_2y)-b^S(x))dy.
\end{split}
\end{equation*}
Hence
\begin{equation}\label{eq:5.7}
\begin{split}
|\partial_x^\beta b_{\theta_2}(x)|&\le\theta_2^{-1}\int_{\mathbb{R}^d}|(\partial_y^\beta\phi)(y)||b^S(x+\theta_2y)-b^S(x)|dy\\
&\lesssim\theta_2^{-1+\alpha'}(r+1)^{-3/2-\rho_1}\int_{\mathbb{R}^d} |y|^{\alpha'}|(\partial_y^\beta\phi)(y)|dy\\
&\lesssim\theta_2^{-1+\alpha'}(r+1)^{-3/2-\rho_1}.
\end{split}
\end{equation}
By (\ref{eq:5.5}) and (\ref{eq:5.7}) we deduce that
\begin{equation}\label{eq:5.8}
|\partial_r b_{\theta_2}(rw)|
\lesssim\theta_2^{-1+\alpha'}(r+1)^{-3/2-\rho_1},
\end{equation}
\begin{equation}\label{eq:5.9}
|\partial_x^\beta(x\cdot b_{\theta_2}(x))|
\lesssim\theta_2^{-1+\alpha'}(r+1)^{-1/2-\rho_1}+(r+1)^{-\rho_2}.
\end{equation}
It is also easy to see that (\ref{eq:1.6}) implies
\begin{equation}\label{eq:5.10}
|{\rm div}\,(b^S(x)-b_{\theta_2}(x))|\lesssim (r+1)^{-1-\rho}.
\end{equation}
We will make use of the Carleman estimates obtained in Section 3 under the condition 
(\ref{eq:1.8}) with $V_S$, $V_L$, $b^S$ and $b^L$ replaced by
$V_S-V_{\theta_1}$, $V_L+V_{\theta_1}$, $b^S-b_{\theta_2}$ and $b^L+b_{\theta_2}$, respectively.
Using the estimates (\ref{eq:5.1}), (\ref{eq:5.3})--(\ref{eq:5.6}) and  (\ref{eq:5.8})--(\ref{eq:5.10}) one can easily check that our assumptions guarantee that 
the new functions $\widetilde V_S$, $\widetilde V_L$ and $\sigma$ satisfy the bounds
\begin{equation*}
\begin{split}
|\widetilde V_S|&\lesssim\left(\theta_1^{\alpha}+\theta_2^{\alpha'}\right)(r+1)^{-3-\rho}+h(r+1)^{-1-\rho},\\
\partial_r(\widetilde V_L-\sigma^2)&\lesssim (r+1)^{-1-\rho}+\left(\theta_1^{-1+\alpha}+\theta_2^{-1+\alpha'}\right)(r+1)^{-3-\rho},
\end{split}
\end{equation*}
with some $\rho>0$. 
Moreover, in view of (\ref{eq:5.2}), $\widetilde V_L$ satisfies (\ref{eq:1.1}) with a new decreasing function $p$ independent of the
parameters $\theta_1$ and $\theta_2$. Then the functions $B_j$ in this case become
\begin{equation*}
\begin{split}
B_1(r)=&\left(\theta_1^{2\alpha}+\theta_2^{2\alpha'}\right)(r+1)^{-6-\rho}+h^2\theta_2^{2\alpha'}r^{-2}(r+1)^{-3-\rho}\\
&+ \theta_2^{2\alpha'}(r+1)^{-3-\rho}\varphi'^2+h^2\varphi''^2+h^2(r+1)^{-2-\rho},\\
B_2(r)=& \theta_2^{2\alpha'}(r+1)^{-3-\rho},\\
B_3(r)=& (r+1)^{-\rho},\\
B_4(r)=&(r+1)^{-1-\rho}+\theta_2^{-1+\alpha'}(r+1)^{-3/2-\rho},\\
B_5(r)=&(r+1)^{-\rho}+\theta_2^{-1+\alpha'}(r+1)^{-2-\rho}.
\end{split}
\end{equation*}
For the new functions $A_1$ and $A_2$ we obtain
\begin{equation*}
\begin{split}
A_1(r)&=\tau^2(\psi'^2\mu)'+(E-\widetilde V_L)\mu'-C\mu r^{-1}(r+1)^{-\rho}\\
&-C\left(M_2\theta_1^{-1+\alpha}+M_2\theta_2^{-1+\alpha'}+M_1\theta_2^{-2+2\alpha'}\right)(r+1)^{-1-\rho}\\
&-Ch^{-1}\tau^{-1}\left(\theta_1^{2\alpha}+\theta_2^{2\alpha'}\right)(r+1)^{-6-\rho}
\mu(\psi')^{-1}-C\tau h^{-1}\theta_2^{2\alpha'}(r+1)^{-3-\rho}\mu\psi'\\
&-Ch\tau^{-1}(r+1)^{-2-\rho}\mu(\psi')^{-1}-Ch\tau\mu(\psi'')^2(\psi')^{-1}\\
&-C\theta_2^{2\alpha'}\mu^2(r\mu')^{-2}(r+1)^{-3-\rho}\mu'-C(rB_3)^2M_1-ChB_3M_2,\\
2A_2(r)&=2r^{-1}\mu-\mu'-Ch^{-1}\tau^{-1}\theta_2^{2\alpha'}(r+1)^{-3-\rho}\mu(\psi')^{-1},
\end{split}
\end{equation*}
where the functions $M_1$ and $M_2$ are the same as in Section 3.
 Note that for $1<r<a$ we have 
 \begin{equation*}
 \mu(\psi')^{-1}=\mathcal{O}(r^3),\quad \mu\psi'=\mathcal{O}(r),\quad \mu(\psi'')^2(\psi')^{-1}=\mathcal{O}(r^{-1}),
 \quad \mu^2(r\mu')^{-2}=\mathcal{O}(1).
 \end{equation*}
Set
\begin{equation*}
\chi_1=\tau^{-2}\theta_1^{-1+\alpha},\quad\eta_1=h^{-1}\tau^{-3}\theta_1^{2\alpha},\quad 
\chi_2=\tau^{-2}\theta_2^{-2+2\alpha'},\quad\eta_2=h^{-1}\tau^{-1}\theta_2^{2\alpha'}.
\end{equation*}
We take $\theta_1=h^{2/(3+\alpha)}$, $\theta_2=h^{1/(1+\alpha')}$, $\tau=h^{-n+1}\tau_0$, where
$\tau_0\gg 1$ is independent of $h$. It is easy to check that with this choice we have the inequalities
\begin{equation*}
\chi_1\le\tau_0^{-2},\quad\eta_1\le\tau_0^{-3},\quad 
\chi_2\le\tau_0^{-2},\quad\eta_2\le\tau_0^{-1}.
\end{equation*}
Therefore, we can make all these parameters small enough by taking $\tau_0$ big enough. 
On the other hand, if we choose $\nu<\rho$, $\ell<\rho$ and 
the above parameters are small enough, it is easy to see that we have $A_1\ge E\mu'/3$
and $A_2\ge 0$. Indeed, in view of (\ref{eq:3.14}), for $r<a$, $r\neq 1$, we have 
\begin{equation*}
\begin{split}
A_1(r)&\ge \tau^2(C_\nu-C(\chi_1+\eta_1+\chi_2+\eta_2))(r+1)^{-1-\nu}\\
&+(E-\widetilde V_L)\mu'-C(r+1)^{-\rho}\mu'-C(h\tau+\theta_2^{2\alpha'})\mu'\\
 &\ge 2^{-1}\tau^2C_\nu(r+1)^{-1-\nu}+(E-\widetilde V_L)\mu'-C(r+1)^{-\rho}\mu'-C(h\tau+\theta_2^{2\alpha'})\mu',\\
 A_2(r)&\ge C_\ell-C\eta_2>0.
 \end{split}
\end{equation*}
 Since $n<2$, we can arrange 
 \begin{equation*}
 C(h\tau+\theta_2^{2\alpha'})\lesssim h^{2-n}+h^{2\alpha'/(1+\alpha')}\le E/4
 \end{equation*} 
 by taking $h$ small enough. Therefore, by (\ref{eq:3.15}) and (\ref{eq:3.16})
 we obtain the desired lower bound for $A_1$ in this case. For $r>a$ the lower bounds for $A_1$ and $A_2$ are obtained in the same
 way as in Section 3 taking $a=h^{-m}$ with $m>0$ big enough. Thus we conclude that the estimate (\ref{eq:3.19}) still holds 
in this case with the new parameters $a,\tau$, as desired.
 \eproof
 
 \section{Carleman estimates for $L^\infty$ potentials $V_S$ and H\"older potentials $b^S$}

In this section we prove the following

\begin{Theorem} \label{6.1}
Under the conditions of Theorem \ref{1.3}, the estimate (\ref{eq:3.1}) holds true with new parameters $a=h^{-m}$, 
$\tau=\tau_0h^{-\frac{1-\alpha'}{1+\alpha'}}$ if $\delta>(1-\alpha')^{-1}$, $\alpha'\le\frac{1}{2}$, and 
$\tau=h^{-\frac{1}{2\delta-1}-\epsilon}\tau_0$, 
 $0<\epsilon\ll 1$ being arbitrary, independent of $h$, if $\alpha'\ge 1-\delta^{-1}$, $1<\delta\le 2$, 
with some positive constants $m$ and $\tau_0$.
\end{Theorem}

{\it Proof.} In this case we will make use of the Carleman estimates obtained in Section 3 under the condition 
(\ref{eq:1.8}) with $b^S$ and $b^L$ replaced by $b^S-b_{\theta_2}$ and $b^L+b_{\theta_2}$, respectively.
The functions $V_S$ and $V_L$ remain unchanged. 
Then the new functions $\widetilde V_S$, $\widetilde V_L$ and $\sigma$ satisfy the bounds
\begin{equation*}
\begin{split}
&|\widetilde V_S|\lesssim (r+1)^{-\delta}+\theta_2^{\alpha'}(r+1)^{-3-\rho}+h(r+1)^{-1-\rho},\\
&\partial_r(\widetilde V_L-\sigma^2)\lesssim (r+1)^{-1-\rho}+\theta_2^{-1+\alpha'}(r+1)^{-3-\rho}.
\end{split}
\end{equation*}
Moreover, $\widetilde V_L$ satisfies (\ref{eq:1.1}) with a new decreasing function $p$ independent of the
parameter $\theta_2$. Then the functions $B_j$ in this case become
\begin{equation*}
\begin{split}
B_1(r)=&(r+1)^{-2\delta}+\theta_2^{2\alpha'}(r+1)^{-6-\rho}+h^2\theta_2^{2\alpha'}r^{-2}(r+1)^{-3-\rho}\\
&+ \theta_2^{2\alpha'}(r+1)^{-3-\rho}\varphi'^2+h^2\varphi''^2+h^2(r+1)^{-2-\rho},\\
B_2(r)=&\theta_2^{2\alpha'}(r+1)^{-3-\rho},\\
B_3(r)=&(r+1)^{-\rho},\\
B_4(r)=&(r+1)^{-1-\rho}+\theta_2^{-1+\alpha'}(r+1)^{-3/2-\rho},\\
B_5(r)=&(r+1)^{-\rho}+\theta_2^{-1+\alpha'}(r+1)^{-2-\rho}.
\end{split}
\end{equation*}
For the new functions $A_1$ and $A_2$ we obtain
\begin{equation*}
\begin{split}
A_1(r)&=\tau^2(\psi'^2\mu)'+(E-\widetilde V_L)\mu'-C\mu r^{-1}(r+1)^{-\rho}\\
&-C\left(M_2\theta_2^{-1+\alpha'}+M_1\theta_2^{-2+2\alpha'}\right)(r+1)^{-1-\rho}\\
&-Ch^{-1}\tau^{-1}(r+1)^{-2\delta}\mu(\psi')^{-1}-Ch^{-1}\tau^{-1}\theta_2^{2\alpha'}(r+1)^{-6-\rho}
\mu(\psi')^{-1}\\
&-C\tau h^{-1}\theta_2^{2\alpha'}(r+1)^{-3-\rho}\mu\psi'-Ch\tau^{-1}(r+1)^{-2-\rho}\mu(\psi')^{-1}-Ch\tau\mu(\psi'')^2(\psi')^{-1}\\
&-C\theta_2^{2\alpha'}\mu^2(r\mu')^{-2}(r+1)^{-3-\rho}\mu'-C(rB_3)^2M_1-ChB_3M_2,\\
2A_2(r)&=2r^{-1}\mu-\mu'-Ch^{-1}\tau^{-1}\theta_2^{2\alpha'}(r+1)^{-3-\rho}\mu(\psi')^{-1}.
\end{split}
\end{equation*}
We have to show that the functions $A_1$ and $A_2$ satisfy the same lower bounds as in the previous section. It suffices to do so
for $r<a$ since for $r>a$ the lower bounds for $A_1$ and $A_2$ can be obtained in the same
 way as in Section 3 taking $a=h^{-m}$ with $m>0$ big enough. 
Set
\begin{equation*}
 \chi=h^{-1}\tau^{-3},\quad 
\chi_2=\tau^{-2}\theta_2^{-2+2\alpha'},\quad\eta_2=h^{-1}\tau^{-1}\theta_2^{2\alpha'}.
\end{equation*}
We take $\theta_2=h^{\frac{1}{1+\alpha'}}$, $\tau=h^{-\frac{1-\alpha'}{1+\alpha'}}\tau_0$, where
$\tau_0\gg 1$ is independent of $h$. Then we have 
\begin{equation*}
\chi_2=\tau_0^{-2},\quad\eta_2=\tau_0^{-1}.
\end{equation*}
Consider first the case $\alpha'=\frac{1}{2}$ and $\delta>2$. Then  
\begin{equation*}
\chi=\tau_0^{-3}.
\end{equation*}
Choose $\nu<\min\{\rho,2(\delta-2)\}$. 
For $r<a$, $r\neq 1$, we have 
\begin{equation*}
\begin{split}
A_1(r)&\ge \tau^2(C_\nu-C(\chi+\chi_2+\eta_2))(r+1)^{-1-\nu}\\
&+(E-\widetilde V_L)\mu'-C(r+1)^{-\rho}\mu'-C(h\tau+\theta_2^{2\alpha'})\mu'\\
 &\ge 2^{-1}\tau^2C_\nu(r+1)^{-1-\nu}+(E-\widetilde V_L)\mu'-C(r+1)^{-\rho}\mu'-C(h\tau+\theta_2^{2\alpha'})\mu',\\
 A_2(r)&\ge C_\ell-C\eta_2>0,
 \end{split}
\end{equation*}
 if $\tau_0$ is taken big enough. In the same way as in the previous section one can conclude that $A_1\ge E\mu'/3$, 
 which leads to the desired estimate in this case.
 
 Let now $\alpha'<\frac{1}{2}$ and $(1-\alpha')^{-1}<\delta\le 2$. Let $1<T<a$ be a parameter of the form 
 \begin{equation*}
 T=T_0(h\tau)^{-\frac{1}{2(\delta-1)}},
 \end{equation*}
  where $T_0\gg 1$ is independent of $h$ and $\tau_0$. For $T\le r<a$ we have 
 \begin{equation*}
 Ch^{-1}\tau^{-1}(r+1)^{-2\delta}\mu(\psi')^{-1}\lesssim h^{-1}\tau^{-1}(r+1)^{-2\delta+2}\mu'(r)\lesssim T_0^{-2\delta+2}\mu'(r)\le
 4^{-1}E\mu'(r),
 \end{equation*}
 provided $T_0$ is taken big enough. On the other hand, for $r<T$, we have
 \begin{equation*}
\begin{split}
 &Ch^{-1}\tau^{-1}(r+1)^{-2\delta}\mu(\psi')^{-1}\lesssim h^{-1}\tau^{-1}(r+1)^{-2\delta+3}\\
 &\lesssim 
  \chi\tau^2(r+1)^{-1-\nu}T^{4-2\delta+\nu}\lesssim 
  \widetilde\chi\tau^2(r+1)^{-1-\nu},
  \end{split}
\end{equation*}
  where
  \begin{equation*}
  \widetilde\chi=\chi(h\tau)^{-\frac{2-\delta+\nu/2}{\delta-1}}=
  h^{-\frac{1+\nu/2}{\delta-1}}\tau^{-\frac{2\delta-1+\nu/2}{\delta-1}}\le \tau_0^{-\frac{2\delta-1+\nu/2}{\delta-1}},
  \end{equation*}
 provided $\nu$ is chosen small enough. Hence we can make $\widetilde\chi$ small enough by taking $\tau_0$ big enough, depending on
 $T_0$. Thus we conclude that the above lower bound for $A_1$ holds with $\chi$ replaced by $\widetilde\chi$
 and $E$ replaced by $3E/4$, which again leads to the inequality $A_1\ge E\mu'/3$ in this case. 
 
 Let now $\alpha'\ge 1-\delta^{-1}$, $1<\delta\le 2$. In this case we take $\tau=h^{-\frac{1}{2\delta-1}-\epsilon}\tau_0$, where
 $0<\epsilon\ll 1$ is arbitrary, independent of $h$, and 
$\tau_0\gg 1$ is independent of $h$ and $\epsilon$. Let $\theta_2$ be as above. Clearly, 
we have $\tau\ge h^{-\frac{1-\alpha'}{1+\alpha'}}\tau_0$,
which implies the inequalities
 \begin{equation*}
\chi_2\le\tau_0^{-2},\quad\eta_2\le\tau_0^{-1}.
\end{equation*}
Moreover, it is easy to see that the above bound of $\widetilde\chi$ still holds, provided $\nu$ is chosen small enough
depending on $\epsilon$. Therefore, the inequality $A_1\ge E\mu'/3$ holds in this case, too.
 
 \eproof

 \section{Resolvent estimates}
 
In this section we show that the Carleman estimates (\ref{eq:3.1}), proved under various conditions,
  imply the resolvent bounds in Theorems \ref{1.1}, \ref{1.2} and \ref{1.3}. The key points are the symmetry of the
  operator $P(h)$ on the Hilbert space $L^2(\mathbb{R}^d)$ and the fact that the potentials belong to $L^\infty$. 
  We need the following lemma the proof of which is given in the next section.

\begin{lemma} \label{7.1} 
Let $\alpha$ and $\beta$ be multi-indices such that $|\alpha|\le 1$, $|\beta|\le 1$, and let $k\in\{0,1\}$ be such that
$k+|\alpha|+|\beta|\le 2$. If $z\in\mathbb{C}$, ${\rm Im}\,z\neq 0$, then the operator
\begin{equation*}
\partial_x^\alpha(\zeta\partial_r)^k(P(h)-z)^{-1}\partial_x^\beta: L^2\to L^2
\end{equation*}
is bounded.
\end{lemma}

Given any $f_1\in L^2$, set
\begin{equation*}
f=(P(h)-E\pm i\varepsilon)^{-1}\langle x\rangle^{-s}f_1\in D(P(h)).
 \end{equation*} 
 By Lemma \ref{7.1} with $\beta=0$, we have $(\zeta\partial_r)^kf\in H^1$, $k=0,1$. We now apply Theorem \ref{3.1} to the function $f$. 
 By (\ref{eq:3.1}) and Lemma \ref{2.2} we have the estimate 
\begin{equation}\label{eq:7.1}
\|\langle x\rangle^{-s}f\|_{L^2}\le N_1\|\langle x\rangle^{s}(P(h)-E\pm i\varepsilon)f\|_{L^2}
+N_2\varepsilon^{1/2}\left(\|f\|_{L^2}+\tau^{-1}\|h\nabla f\|_{L^2}\right)
\end{equation}
where 
 \begin{equation*}
N_1=ah^{-1}\exp\left(Ch^{-1}\tau\log a\right),\quad N_2=(h\tau)^{1/2}N_1,
 \end{equation*}
with a constant $C>0$ independent of $h$ and $\varepsilon$. On the other hand, since the operator $P(h)$ is symmetric, we have
\begin{equation*}
\begin{split}
\varepsilon\|f\|^2_{L^2}&=\pm{\rm Im}\,\langle (P(h)-E\pm i\varepsilon)f,f\rangle_{L^2}\\
&\le \gamma^2N_2^{-2}\|\langle x\rangle^{-s}f\|^2_{L^2}+\gamma^{-2}N_2^2\|\langle x\rangle^{s}(P(h)-E\pm i\varepsilon)f\|^2_{L^2}
\end{split}
\end{equation*}
for every $\gamma>0$, which yields
\begin{equation}\label{eq:7.2}
N_2\varepsilon^{1/2}\|f\|_{L^2}\le \gamma\|\langle x\rangle^{-s}f\|_{L^2}+
\gamma^{-1}N_2^2\|\langle x\rangle^{s}(P(h)-E\pm i\varepsilon)f\|_{L^2}.
\end{equation} 
We also have
\begin{equation*}
\begin{split}
{\rm Re}\,\langle (P(h)-E\pm i\varepsilon)f,f\rangle_{L^2}&=\|(ih\nabla+b)f\|^2_{L^2}+{\rm Re}\,\langle (V-E)f,f\rangle_{L^2}\\
&\ge \|h\nabla f\|^2_{L^2}-C\|f\|^2_{L^2}.
\end{split}
\end{equation*}
Thus we get
\begin{equation}\label{eq:7.3}
\|h\nabla f\|_{L^2}\lesssim \|f\|_{L^2}+\|(P(h)-E\pm i\varepsilon)f\|_{L^2}.
\end{equation}
Taking $\gamma$ small enough, independent of $N_2$ and $\varepsilon$, by (\ref{eq:7.1}), (\ref{eq:7.2}) and (\ref{eq:7.3}) we obtain
the estimate
\begin{equation}\label{eq:7.4}
\|\langle x\rangle^{-s}f\|_{L^2}\lesssim (N_1+N_2^2)\|\langle x\rangle^{s}(P(h)-E\pm i\varepsilon)f\|_{L^2}.
\end{equation}
It follows from (\ref{eq:7.4}) that the resolvent estimate
\begin{equation}\label{eq:7.5}
\left\|\langle x\rangle^{-s}(P(h)-E\pm i\varepsilon)^{-1}\langle x\rangle^{-s}
\right\|_{L^2\to L^2}\lesssim N_1+h\tau N_1^2\lesssim N_1^2
\end{equation}
holds for all $0<h\ll 1$. Clearly,
(\ref{eq:7.5}) implies the desired bounds for $g_s^\pm$.

\section{Proof of Lemma \ref{7.1}}

It follows from the resolvent identity
\begin{equation*}
(P(h)-z)^{-1}=(P(h)-z_0)^{-1}+(z-z_0)(P(h)-z_0)^{-1}(P(h)-z)^{-1}
\end{equation*}
that if the lemma holds for one $z_0\in\mathbb{C}$ with ${\rm Im}\,z_0\neq 0$, then it holds for all
$z\in\mathbb{C}$ with ${\rm Im}\,z\neq 0$. We will prove it for $z_0=\pm i\theta$ with real $\theta\gg 1$. 
Set $P_0(h)=-h^2\Delta$. 
The self-adjoint realization of the operator $-h^2\Delta$ on $L^2(\mathbb{R}^d)$ will be again denoted by $P_0(h)$.
Clearly, the lemma holds with $P(h)$ replaced by $P_0(h)$. 

We first prove the lemma with $k=0$. We have
\begin{equation*}
P(h)=P_0(h)+ihb\cdot\nabla+ih\nabla\cdot b+|b|^2+V.
\end{equation*}
Using the resolvent identity
\begin{equation*}
(P(h)\pm i\theta)^{-1}=(P_0(h)\pm i\theta)^{-1}-(P_0(h)\pm i\theta)^{-1}(P(h)-P_0(h))(P(h)\pm i\theta)^{-1}
\end{equation*}
we get
\begin{equation}\label{eq:8.1}
(P(h)\pm i\theta)^{-1}=\mathcal{T}_1+\mathcal{T}_2ihb\cdot\nabla(P(h)\pm i\theta)^{-1},
 \end{equation}
 where
 \begin{align*}
 &\mathcal{T}_1=(P_0(h)\pm i\theta)^{-1}-(P_0(h)\pm i\theta)^{-1}(ih\nabla\cdot b+|b|^2+V)(P(h)\pm i\theta)^{-1},\\
 &\mathcal{T}_2=-(P_0(h)\pm i\theta)^{-1}.
 \end{align*}
 By (\ref{eq:8.1}) we have
 \begin{equation}\label{eq:8.2}
ihb\cdot\nabla(P(h)\pm i\theta)^{-1}=\mathcal{K}_1+\mathcal{K}_2ihb\cdot\nabla(P(h)\pm i\theta)^{-1},
 \end{equation}
 where
 \begin{equation*}
 \mathcal{K}_j=ihb\cdot\nabla\mathcal{T}_j,\quad j=1,2.
 \end{equation*}
 Observe now that
 \begin{equation*}
 \begin{split}
 \|\mathcal{K}_2\|_{L^2\to L^2}&\lesssim \|P_0(h)^{1/2}(P_0(h)\pm i\theta)^{-1}\|_{L^2\to L^2}\\
 &\lesssim\sup_{\sigma\ge 0}|\sigma(\sigma^2\pm i\theta)^{-1}|\lesssim\theta^{-1/2}\le 1/2
 \end{split}
 \end{equation*}
 if $\theta$ is big enough. Therefore, the operator $I-\mathcal{K}_2$ is invertible on $L^2$, so the identity (\ref{eq:8.2})
 can be written in the form
 \begin{equation}\label{eq:8.3}
ihb\cdot\nabla(P(h)\pm i\theta)^{-1}=(I-\mathcal{K}_2)^{-1}\mathcal{K}_1.
 \end{equation}
 By (\ref{eq:8.1}) and (\ref{eq:8.3}) we obtain
 \begin{equation}\label{eq:8.4}
(P(h)\pm i\theta)^{-1}=\mathcal{T}_1+\mathcal{T}_2(I-\mathcal{K}_2)^{-1}\mathcal{K}_1.
 \end{equation}
 Let $|\alpha|=|\beta|=1$. Clearly, the operators $\partial_x^\alpha\mathcal{T}_j:L^2\to L^2$, $j=1,2$, are bounded. Therefore, by 
 (\ref{eq:8.4}) the operator $\partial_x^\alpha(P(h)\pm i\theta)^{-1}:L^2\to L^2$ is bounded, too.
 This implies that the operator 
 \begin{equation*}
 (P(h)\pm i\theta)^{-1}\partial_x^\beta=-(\partial_x^\beta(P(h)\mp i\theta)^{-1})^*:L^2\to L^2
 \end{equation*}
 is bounded. This in turn implies the boundedness of the operators $\mathcal{K}_1\partial_x^\beta:L^2\to L^2$
 and $\partial_x^\alpha\mathcal{T}_1\partial_x^\beta:L^2\to L^2$.
 Thus, by (\ref{eq:8.4}) we conclude that the operator $\partial_x^\alpha(P(h)\pm i\theta)^{-1}\partial_x^\beta:L^2\to L^2$ is bounded,
 as desired.
 
Let now $k=1$.
If $\alpha=0$, then the lemma follows from the previous case. In what follows we will consider the case when $|\alpha|=1$. Then
$\beta=0$. To prove the lemma in this case we are going to use the identity (\ref{eq:8.4}) again. Observe first that we can write
\begin{equation*}
\nabla\cdot b=\nabla\cdot b^L+b^S\cdot\nabla+{\rm div}\,b^S.
\end{equation*}
Therefore, since the operators
\begin{equation*}
\zeta\partial_r(P_0(h)-i\theta)^{-1},\quad (P(h)-i\theta)^{-1}: L^2\to H^1
\end{equation*}
are bounded, 
in view of (\ref{eq:8.4}), it suffices to show that 
the operator
\begin{equation*}
\mathcal{T}_3:=\zeta\partial_r(P_0(h)-i\theta)^{-1}\nabla\cdot b^L(P(h)-i\theta)^{-1}: L^2\to H^1
\end{equation*}
is bounded. To this end we will use the identity
\begin{equation}\label{eq:8.5}
\zeta\partial_r(P_0(h)-i\theta)^{-1}=-(P_0(h)-i\theta)^{-1}[\zeta\partial_r,P_0(h)](P_0(h)-i\theta)^{-1}+(P_0(h)-i\theta)^{-1}\zeta\partial_r.
 \end{equation}
 By (\ref{eq:8.5}) we can write $\mathcal{T}_3=\mathcal{T}_4+\mathcal{T}_5$, where
 \begin{align*}
&\mathcal{T}_4=-(P_0(h)-i\theta)^{-1}[\zeta\partial_r,P_0(h)](P_0(h)-i\theta)^{-1}\nabla\cdot b^L(P(h)-i\theta)^{-1},\\
&\mathcal{T}_5=(P_0(h)-i\theta)^{-1}\zeta\partial_r\nabla\cdot b^L(P(h)-i\theta)^{-1}.
\end{align*}
 Let us now see that the commutator $[\zeta\partial_r,P_0(h)]$ can be written in the form
\begin{equation}\label{eq:8.6}
[\zeta\partial_r,P_0(h)]=\sum_{|\alpha|\le 1}\sum_{|\beta|\le 1}\partial_x^\alpha a_{\alpha,\beta}(x)\partial_x^\beta
\end{equation}
with bounded functions $a_{\alpha,\beta}$. We have
\begin{equation*}
\zeta\partial_r=\sum_{\nu=1}^d\omega_\nu(x)\partial_{x_\nu},\quad \omega_\nu=x_\nu(1+|x|^2)^{-1/2},
\end{equation*}
\begin{equation*}
[\zeta\partial_r,\partial_{x_j}]=-\sum_{\nu=1}^d\frac{\partial \omega_\nu}{\partial x_j}\partial_{x_\nu},
\end{equation*}
\begin{equation*}
\begin{split}
[\zeta\partial_r,\partial_{x_j}^2] &=[\zeta\partial_r,\partial_{x_j}]\partial_{x_j}+\partial_{x_j}[\zeta\partial_r,\partial_{x_j}]\\
 &=-\sum_{\nu=1}^d\frac{\partial \omega_\nu}{\partial x_j}\partial_{x_\nu}\partial_{x_j}
 -\partial_{x_j}\sum_{\nu=1}^d\frac{\partial \omega_\nu}{\partial x_j}\partial_{x_\nu}\\
 &=-\sum_{\nu=1}^d\partial_{x_\nu}\frac{\partial \omega_\nu}{\partial x_j}\partial_{x_j}+
 \sum_{\nu=1}^d\frac{\partial^2 \omega_\nu}{\partial x_\nu\partial x_j}\partial_{x_j}
 -\sum_{\nu=1}^d\partial_{x_j}\frac{\partial \omega_\nu}{\partial x_j}\partial_{x_\nu},
\end{split}
\end{equation*}
which confirms (\ref{eq:8.6}). It follows from the lemma with $k=0$ and (\ref{eq:8.6}) that the operator $\mathcal{T}_4: L^2\to H^1$ is bounded. 
 To prove this for the 
operator $\mathcal{T}_5$ it suffices to show that the operator $\zeta\partial_r\nabla\cdot b^L$ can be written as the right-hand side of 
(\ref{eq:8.6}). To do so it suffices to show that 
the commutator $[\zeta\partial_r,\nabla\cdot b^L]$ can be written in the form
\begin{equation}\label{eq:8.7}
[\zeta\partial_r,\nabla\cdot b^L]=\sum_{|\alpha|\le 1}\partial_x^\alpha a_{\alpha}(x)
\end{equation}
with bounded functions $a_{\alpha}$. We have 
\begin{equation*}
\begin{split}
[\zeta\partial_r,\partial_{x_j}b_j^L] &=\partial_{x_j}\zeta\partial_rb_j^L
+[\zeta\partial_r,\partial_{x_j}]b_j^L-\partial_{x_j}b_j^L\zeta\partial_r\\
&=\partial_{x_j}\left(b_j^L\zeta\partial_r+\zeta\frac{\partial b_j^L}{\partial r}\right)
+[\zeta\partial_r,\partial_{x_j}]b_j^L-\partial_{x_j}b_j^L\zeta\partial_r\\
&=\partial_{x_j}\zeta\frac{\partial b_j^L}{\partial r}
+[\zeta\partial_r,\partial_{x_j}]b_j^L\\
&=\partial_{x_j}\zeta\frac{\partial b_j^L}{\partial r}
-\sum_{\nu=1}^d\partial_{x_\nu}\frac{\partial \omega_\nu}{\partial x_j}b_j^L
+\sum_{\nu=1}^d\frac{\partial^2 \omega_\nu}{\partial x_\nu\partial x_j}b_j^L.
\end{split}
\end{equation*}
Since the functions $\partial_r^\ell b_j^L$, $\ell=0,1$, are bounded by assumption, the above formula confirms (\ref{eq:8.7}).
\eproof

\begin{rem} \label{8.1}
The above proof does not work with $\zeta$ replaced by $1$ because the singularity of the derivatives of the coefficients of the
operator $\partial_r=\sum_{j=1}^dw_j(x)\partial_{x_j}$ at $x=0$ is too strong. 
\end{rem}

\begin{rem} \label{8.2}
In the same way as above, we can prove that the lemma holds with $k=0$ and all multi-indices $\alpha$ and $\beta$
such that $|\alpha|+|\beta|\le 2$, provided we suppose in addition that the function ${\rm div}\,b^L$ exists
and is bounded. However, it is clear from the proof of Theorem \ref{3.1} that we do not need so much regularity in our analysis.
Therefore, we do not need to impose this condition on $b^L$. 
\end{rem}

\appendix
\section{Proof of Lemma \ref{3.2}}

Differentiating the identity $r^2=\sum_{j=1}^dx_j^2$ yields
 \begin{equation*}
 2r\frac{\partial r}{\partial x_j}=2x_j,
 \end{equation*}
 which implies 
 \begin{equation}\label{eq:A.1}
 \frac{\partial r}{\partial x_j}=w_j,\quad 1\le j\le d.
 \end{equation}
 Let $U\subset\mathbb{R}^{d-1}$ be a small open domain and let $Y:U\to \mathbb{S}^{d-1}$ be a diffeomorphisme. 
 Let $y\in U$ be local coordinates and let $w=Y(y)\in Y(U)$. Then we have
 \begin{equation}\label{eq:A.2}
 \partial_{x_j}=\frac{\partial r}{\partial x_j}\partial_r+\frac{\partial y}{\partial x_j}\cdot\partial_y=
 w_j\partial_r+r^{-1}\frac{\partial y}{\partial w_j}\cdot\partial_y.
 \end{equation}
 It is easy to see that (\ref{eq:3.3}) follows from (\ref{eq:A.1}) and (\ref{eq:A.2}). To see that the operator $q_j$
 is antisymmetric with respect to the scalar product
in $L^2(\mathbb{S}^{d-1})$, we will use that the operator $\partial_{x_j}$
 is antisymmetric with respect to the scalar product
in $L^2(\mathbb{R}^{d})=L^2(\mathbb{R}^+\times\mathbb{S}^{d-1},r^{d-1}drdw)$. 
Indeed, this implies that the operator $r^{(d-1)/2}\partial_{x_j} r^{-(d-1)/2}$, acting on functions
$u\in H^1(\mathbb{R}^+\times\mathbb{S}^{d-1},drdw)$, $u(0)=0$, 
 is antisymmetric with respect to the scalar product
in $L^2(\mathbb{R}^+\times\mathbb{S}^{d-1},drdw)$. Therefore, we have the identity
\begin{equation}\label{eq:A.3}
\begin{split}
&\int_0^\infty\int_{\mathbb{S}^{d-1}}(w_j\partial_r+r^{-1}q_j)u(r,w)\overline{v(r,w)}drdw\\
&=-\int_0^\infty\int_{\mathbb{S}^{d-1}}u(r,w)(w_j\partial_r+r^{-1}q_j)\overline{v(r,w)}drdw
\end{split}
\end{equation}
for all $u,v\in H^1(\mathbb{R}^+\times\mathbb{S}^{d-1},drdw)$, $u(0)=v(0)=0$. We now apply (\ref{eq:A.3})
with $u=u_1(r)u_2(w)$, $v=v_1(r)v_2(w)$, $u_1, v_1\in H^1(\mathbb{R}^+)$, $u_1(0)=v_1(0)=0$, and 
$u_2, v_2\in H^1(\mathbb{S}^{d-1})$. Since
\begin{equation*}
\int_0^\infty\partial_ru_1(r)\overline{v_1(r)}dr=-\int_0^\infty u_1(r)\partial_r\overline{v_1(r)}dr,
\end{equation*}
we deduce from (\ref{eq:A.3}),
\begin{equation*}
\int_{\mathbb{S}^{d-1}}q_ju_2(w)\overline{v_2(w)}dw=-\int_{\mathbb{S}^{d-1}}u_2(w)q_j\overline{v_2(w)}dw
\end{equation*}
for all $u_2,v_2\in H^1(\mathbb{S}^{d-1})$, as desired.

\end{document}